\documentclass[10pt]{article}
\usepackage{epsfig,subfigure}
\usepackage{amsmath,amssymb}
\usepackage{exscale}
\usepackage[mathscr]{eucal}
\usepackage{fancyhdr}
\usepackage{graphicx}
\usepackage{multirow,bigstrut}
\usepackage{labelfig}

\usepackage[a4paper]{geometry}
\newtheorem{theorem}{Theorem}

\newtheorem{lemma}[theorem]{Lemma}

\def\RR{{I~\hspace{-1.45ex}R} }

\begin{document}

\begin{center}{\Large \textbf{Second-grade fluids in 
curved pipes}}

\vspace{1cm}
\begin{center}
\textsc{nadir arada}{\footnote{Universit\'e Mohammed Seddik BENYAHIA de Jijel, D\'epartement de Math\'ematiques, Alg\'erie (nadir.arada@gmail.com)}}\hspace{0.5cm}\hspace{0.5cm}
\textsc{paulo correia}{\footnote{Departamento de Matem\'atica and Centro de Investiga\c{c}\~ao em Matem\'atica e Aplica\c{c}\~oes (CIMA), Escola de Ci\^encias e Tecnologia, Universidade de \'Evora, Rua Rom\~ao Ramalho, 59, 7000-651, \'Evora, Portugal (pcorreia@uevora.pt)}} 
\vspace{5mm}
\end{center}

\end{center}

\vspace{1cm}

\begin{abstract}\noindent This paper is concerned with the application of finite element methods to obtain solutions for steady fully developed 
second-grade flows in a curved pipe of circular cross-section and 
arbitrary curvature ratio, under a given axial 
 pressure gradient. The qualitative and quantitative behavior of the secondary flows is analyzed with respect to inertia and viscoelasticity.\vspace{2mm}

\end{abstract}

\markright{\today}

\section{Introduction}
\label{S1} 
Fully developed viscous flow in a curved pipe with circular cross-section was 
first studied theoretically by Dean (\cite{dean1}, \cite{dean2}) applying regular perturbation
methods. He showed that for small curvature
ratio the flow depends only on a single parameter, the so-called {\it Dean number}.
In \cite{sohr}, Soh and Berger  solved the Navier-Stokes equations for the fully developed
flow of an homogeneous Newtonian fluid in a curved pipe of circular cross-section
for arbitrary curvature ratio. They solved numerically the Navier-Stokes system
in the primitive variables form using a finite difference scheme.
Closed form perturbation solutions for a second order model were obtained
by several authors in the special case where the second normal stress coefficient
is zero.\vspace{2mm} 

 For this model, Jitchote and Robertson \cite{citation3} obtained analytical solutions to the
perturbation equations and analyze the effects of non-zero second normal stress 
coefficient on the behaviour of the solution. Theoretical results regarding 
this problem using a splitting method were obtained by Coscia and Robertson 
\cite{CoRo}.\vspace{2mm} 

Our aim here is to apply the finite element method
 to the second-grade model for fully developped flows in a curved pipe  and analyze the non-Newtonian 
effects of the flow. Quantitative and qualitative behaviour of the axial velocity and the stream function
 of creeping and inertial second-grade fluids are also studied.
Similar techniques have been applied to other non-Newtonian models (see \cite{arada-pires-sequeira1}, \cite{arada-pires-sequeira2}). \vspace{2mm}

The paper is organized as follows. In Section $2$, we consider the governing equations and rewrite the problem in an equivalent decoupled form, composed of a Stokes-like system and a transport equation  considered as two auxiliary linear problems. In order to describe the curved pipe geometry, these equations are presented in non-dimensional polar coordinates in Section $3$. In Section $4$, we outline  tha discretization method using a finite element method to obtain approximate solutions to the original problem. Finally, numerical experiments based on a decoupled fixed point algorithm are illustrated in Section $5$. A systematic numerical study of the qualitative behaviour of inertial and viscoelastic effects for a certain range of non-dimensional parameters (Reynolds number, viscoelastic coefficient and pipe curvature) associated to the model is performed. The issue of increasing the values of the parameters for which the convergence of the proposed iterative algorithm can be shown, is addressed by using a continuation method.

\section{Governing equations}
\begin{figure}[!ht]
\centering
\begin{tabular}{c}
\includegraphics[width=42mm]{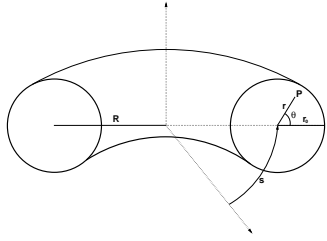}
\end{tabular}
\caption{\small Polar toroidal coordinates in a curved pipe.}
\label{tube}
\end{figure}

We consider steady isothermal flows of incompressible second-grade
fluids in a curved pipe $\Omega\subset \RR^3$ with boundary $\partial \Omega$, 
constant circular cross--section $\Omega_S$ of radius $r$ with the line of centers
 coiled in a circle of radius $R$ (see Figure \ref{tube}).

\noindent
 The corresponding equations are given by
$$
\left\{ \!\!
\begin{array}{ll}
	-\mu \Delta\mathbf{u}-\widetilde\alpha_1\mathbf{u}\cdot\!\nabla\Delta\mathbf{u}
	+\!\rho \mathbf u\cdot\! \nabla\mathbf{u}+\nabla\pi
	\!=\!\nabla\!\cdot \widetilde {\mathbf L}(\mathbf{u})  & \text{ in } \Omega,\vspace{2mm}\\
	\nabla\!\cdot \mathbf{u}=0 & \text{ in } \Omega, \vspace{2mm}\\
	\mathbf{u}=0 & \text{ on } \partial\Omega,
\end{array}
\right.
$$
where $\mathbf{u}$ is the fluid velocity, $\pi$ is the pressure, $\mu$ is the viscosity, $\rho$ is the constant density, $\widetilde\alpha_1$ and 
$\widetilde\alpha_2$ are the Rivlin-Ericksen material constants, and 
$\widetilde{\mathbf{L}}(\mathbf{u})$ is a nonlinear term given by
$$
	\widetilde{\mathbf{L}}(\mathbf{u})\!=\!\widetilde \alpha_1 \nabla\mathbf{u}^{T} 
	\!\left( \nabla\mathbf{u}+\nabla\mathbf{u}^T \right)
	+(\widetilde \alpha_1+\widetilde\alpha_2)\!
	\left(\nabla\mathbf{u}+\nabla\mathbf{u}^T\right)^2.$$
We consider an adimensionalised formulation of the previous system by introducing
the following quantities 
\begin{equation}\label{variadim}
\mathbf x=\textstyle \tfrac{\widetilde {\mathbf x}}{a} ,  \qquad
\mathbf u=\tfrac{\widetilde {\mathbf u}}{U}, \qquad \pi=\tfrac{\widetilde \pi a}{\mu U},\end{equation}
where the symbol $\ \widetilde {}\ $ is attached to dimensional parameters
($U$ represents a characteristic velocity of the flow). We also introduce the Reynolds number $\mathcal{R}e$ and two non--dimensional ratios
involving the constant material moduli $\widetilde\alpha$ and $\widetilde\alpha_2$,
$$
	\mathcal{R}e=\textstyle \tfrac{\rho U a}{\mu}, 
	\qquad \alpha_1=\textstyle\tfrac{U\widetilde\alpha}{\mu a}, 
	\qquad \alpha_2=\textstyle\tfrac{U\widetilde\alpha_2}{\mu a}.$$
The dimensionless system takes the form   
\begin{equation}
\left\{\!\!
\begin{array}[c]{ll}
	-\Delta\mathbf{u}-\alpha_1\mathbf{u}\cdot\!\nabla\Delta\mathbf{u}
	+{\mathcal R}e \, \mathbf u\cdot\! \nabla\mathbf{u}+\nabla\pi
	\!=\!\nabla\!\cdot \mathbf L(\mathbf{u})& \text{ in } \Omega,
	 \vspace{2mm}\\
	\nabla\!\cdot \mathbf{u}=0    & \text{ in } \Omega,
	           \vspace{2mm}\\
	\mathbf{u}=0 & \text{ on } \partial\Omega,
	\end{array}
	\right.  \label{ssfcurved}
\end{equation}
with 
$$\mathbf{L}(\mathbf{u})\!=\!\alpha_1 \nabla\mathbf{u}^{T} 
	\!\big( \nabla\mathbf{u}+\nabla\mathbf{u}^T \big)
	+(\alpha_1+\alpha_2)\!
	\big(\nabla\mathbf{u}+\nabla\mathbf{u}^T\big)^2.$$
System (\ref{ssfcurved}) can be rewritten in the following equivalent
form involving a Stokes system and a transport equation,
\begin{equation}\left\{\!\!
\begin{array}[c]{ll}\label{couplestcurved}
	-\Delta\mathbf{u}+\nabla p
	=\boldsymbol\varrho & \text{ in } \Omega,         \vspace{2mm}\\
         \nabla \cdot \mathbf{u}=0 & \text{ in } \Omega,\vspace{2mm}\\
        \mathbf{u}=0 & \text{ on } \partial\Omega,
	              \vspace{2mm}\\
	\boldsymbol\varrho\!+\!\alpha_1\mathbf{u}\cdot\!\nabla \boldsymbol\varrho
	\!=\!\nabla\!\cdot\!\left(\mathbf L(\mathbf{u})
	\!-\!{\cal R}e \, \mathbf u\otimes \mathbf u\!-\!\alpha_1 p  
	\nabla\!\mathbf{u}^{T}\right) & \text{ in } \Omega,
\end{array}
\right.
\end{equation}
where the new variable $p$ is given by $\pi\!=\!p\!+\!\alpha_1 \mathbf{u}\cdot\!\nabla p$.
In the following lemma, we (formally) prove that if $\mathbf u$ is such that the operator 
$I+\alpha_ 1 \mathbf u\cdot \nabla$ is invertible, then problems (\ref{ssfcurved}) and (\ref{couplestcurved}) are equivalent. 
\begin{lemma} \label{invertibility}
If $(\mathbf u,p,\varrho)$ is a strong solution of $(\ref{couplestcurved})$, then 
$(\mathbf u,p+\alpha \mathbf u\cdot \nabla p)$ is a strong solution of $(\ref{ssfcurved})$. Conversely, if 
$(\mathbf u,\pi)$ is a strong solution for $(\ref{ssfcurved})$ and if 
$I+\alpha_ 1 \mathbf u\cdot \nabla$ is
 invertible, then $(\mathbf u,(I+\alpha_ 1 \mathbf u\cdot \nabla)^{-1}\pi,
\boldsymbol\varrho)$ 
is a solution of $(\ref{couplestcurved})$.
\end{lemma}
{\it Proof.}
Let $(\mathbf{u},p,\boldsymbol\varrho)$ be a solution of 
 (\ref{couplestcurved}). Applying the operator $I+\alpha_ 1 \mathbf u\cdot \nabla$
to both sides of $(\ref{couplestcurved})_1$, and taking into account $(\ref{couplestcurved})_2$-$(\ref{couplestcurved})_4$, we obtain
	$$\begin{array}{ll}
	\big(I+\alpha_ 1 \mathbf u\cdot \nabla\big) \big(\Delta \mathbf u
	+\nabla p\big)=
	\big(I+\alpha_ 1 \mathbf u\cdot \nabla\big)\boldsymbol\varrho
\vspace{2mm}\\
	=\nabla \cdot \big(\mathbf L(\mathbf{u})
	  -{\cal R}e \, \mathbf u\otimes \mathbf u-\alpha_1 p \left(  
	\nabla\mathbf{u}\right)^{T}\big)\vspace{2mm}\\
	=\nabla \cdot \mathbf L(\mathbf{u})-{\cal R}e \,
	 \mathbf u\cdot \nabla\mathbf u
	-\alpha_ 1\nabla \big(\mathbf u\cdot \nabla p\big)
	+\alpha_ 1\nabla \cdot \big(\nabla p\otimes \mathbf u\big)
	\vspace{2mm}\\
	=\nabla \cdot \mathbf L(\mathbf{u})-{\cal R}e \,
	 \mathbf u\cdot \nabla\mathbf u-
	\alpha_ 1\nabla \big(\mathbf u\cdot \nabla p\big)
	+\alpha_ 1\mathbf u\cdot \nabla(\nabla p)\vspace{2mm}\\
	=\nabla \cdot \mathbf L(\mathbf{u})-{\cal R}e \,
	\mathbf u\cdot \nabla\mathbf u-
	\nabla \big(p+\alpha_ 1\mathbf u\cdot \nabla p\big)+
	\big(I+\alpha_ 1\mathbf u\cdot \nabla\big) \nabla p.
	\end{array}$$
Setting $\pi=p +\alpha_1 \mathbf{u}\cdot\nabla p$
 we easily see that $(\mathbf{u},\pi)$ 
is a solution of $(\ref{ssfcurved})$.
Conversely, let $(\mathbf{u},\pi)$ be a solution of (\ref{ssfcurved}). From $(\ref{ssfcurved})_1$, it follows that
	$$
	-\left(I+\alpha_ 1 \mathbf u\cdot \nabla\right) \, \Delta \mathbf u+{\cal R}e \ \mathbf u\cdot \nabla\mathbf u+\nabla \pi=
	\nabla\cdot\mathbf{L}(\mathbf{u}).$$
This equation together with the fact that
	$$\big(I+\alpha_ 1 \mathbf u\cdot \nabla\big)\nabla p-\nabla \pi=-\alpha_1
	\nabla \cdot
	\left(p ( 
	\nabla\mathbf{u})^{T}\right),$$
gives
	$$\big(I+\alpha_ 1 \mathbf u\cdot \nabla\big)\big(-\Delta \mathbf u+
	\nabla p\big)=\nabla\cdot\big(\mathbf{L}(\mathbf{u})-{\cal R}e \, \mathbf u\otimes \mathbf u-\alpha_1p ( 
	\nabla\mathbf{u})^{T}\big).$$
therefore, if the operator $I+\alpha_ 1 \mathbf u\cdot \nabla$ is invertible,
 we deduce that 
	$$-\Delta \mathbf u+
	\nabla p=\big(I+\alpha_ 1 \mathbf u\cdot \nabla\big)^{-1} 
	\big(\nabla\cdot\big(\mathbf{L}(\mathbf{u})-{\cal R}e \, \mathbf u\otimes \mathbf u-\alpha_1
	p ( 
	\nabla\mathbf{u})^{T}\big)\big)\equiv \boldsymbol\varrho.$$
This completes the proof.$\hfill\Box$ \vspace{2mm}\\
Throughout the paper, we consider the  case of $\alpha_1+\alpha_2=0$, corresponding to a thermodynamical compatibility condition. For a more detailed analysis 
of this problem see \cite{correia}. For simplicity, we set $\alpha_1\equiv \alpha$.
\section{Formulation in polar toroidal coordinates}
Since we are interested in studying the behaviour of steady flows in a curved pipe with circular cross-section,
 it is more convenient to use the polar toroidal coordinate system, in the
variables $(\widetilde r,\widetilde\theta,\widetilde s)$,
 defined with respect to the rectangular cartesian coordinates 
$(\widetilde x,\widetilde y,\widetilde z)$ 
through the relations
	$$ \widetilde x=\left(R+\widetilde r\cos\theta\right)
         \cos\tfrac{\widetilde s}{R},\qquad\widetilde y=\left(R+\widetilde r\cos\theta\right) 
	\sin\tfrac{\widetilde s}{R},\qquad
	\widetilde z=\widetilde r\sin\theta,$$
with $0<r_0<R$, $0\leq\widetilde\theta<2\pi$ and $0\leq \widetilde s < \pi R$.
Introducing the axial variable and the pipe curvature ratio
	$$ s=\tfrac{\widetilde s}{r_0}, \qquad \delta=\tfrac{r_0}{R},$$
the  corresponding non-dimensional coordinate system is 
 given by 
	$$
	x=\left(\tfrac{1}{\delta}+r\cos\theta\right) \cos (s\delta),\qquad y=\left(\tfrac{1}{\delta}+r\cos\theta\right) 
	\sin (s\delta),\qquad
	z=r\sin\theta,$$
with $\delta<1$, $0\leq\theta<2\pi$ and $0\leq s < 
\tfrac{\pi}{\delta}$ (Here $\delta$ is the pipe curvature ratio.). Let us now formulate problem (\ref{couplestcurved}) in this
new coordinate system. To simplify the redaction we set 
	$$ B_1=r\delta  \sin\theta,
	\qquad  B_2=r\delta  \cos\theta,\qquad
	 B=1+r\delta \cos\theta.$$
 By using standard arguments, we first rewrite 
the Stokes problem ($\ref{couplestcurved})_{1,2,3}$ in the
toroidal coordinates $(r,\theta,s)$, and obtain
\begin{equation}\left\{
\begin{array}{ll} 
\displaystyle	{\cal A}u+\tfrac{2}{r^{2}}\tfrac{\partial v}{\partial \theta}
	+\tfrac{2B_2}{rB^{2}} \tfrac{\partial w}{\partial s}
	+\tfrac{B^2+B_2^2}{(rB)^2}u-\tfrac{B_1(B+B_2)}{(rB)^2}v
	+\tfrac{\partial p}{\partial r}=\varrho_1,\vspace{2mm}\\
\displaystyle{\cal A}v-\tfrac{2}{r^2}\tfrac{\partial u}{\partial \theta}
	-\tfrac{2B_1}{rB^2}\tfrac{\partial w}{\partial s}
	+\tfrac{B_1}{(rB)^2}u+\tfrac{B^2+B_1^2}{(rB)^2}v+
	\tfrac{1}{r}\tfrac{\partial p}{\partial \theta}=\varrho_2,\vspace{2mm}\\
	\displaystyle{\cal A}w-\tfrac{2B_2}{rB^{2}} 
	\tfrac{\partial u}{\partial s}
	+\tfrac{2B_1}{rB^{2}} \tfrac{\partial v}{\partial s}+\left(\tfrac
	{\delta}{B}\right)^2w+\tfrac{1}{B}\tfrac{\partial p}{\partial s}=\varrho_3,\vspace{2mm}\\
 \displaystyle{\tfrac{\partial u}{\partial r}+\tfrac{1}{r}
	\tfrac{\partial v}{\partial \theta}+\tfrac{1}{B}
	\tfrac{\partial w}{\partial s}+\tfrac{B+B_2}{rB}u-
	\tfrac{B_1}{rB}v=0,}
	  \end{array}\right.\label{stokescurv}
\end{equation}
where 
  $${\cal A}\zeta=-\tfrac{\partial ^2\zeta}{\partial r^2}-\tfrac{1}{r^2}
	\tfrac{\partial ^2\zeta}{\partial \theta^2}-
	\tfrac{1}{B^2}\tfrac{\partial ^2\zeta}{\partial s^2}
	-\tfrac{B+B_2}{rB} 
	\tfrac{\partial \zeta}{\partial r}
	+\tfrac{B_1}{r^2B}
	\tfrac{\partial \zeta}{\partial \theta}.$$
In the same way, the transport equation ($\ref{couplestcurved})_{4}$ 
reads as
\begin{equation}\label{transportcurv}\left\{\begin{array}{lll}
	\displaystyle\varrho_1+\alpha\left(u \,
	 \tfrac{\partial \varrho_1}{\partial r}
	+\tfrac{v}{r} \, \tfrac{\partial \varrho_1}{\partial\theta}
	 +\tfrac{w}{B} \, 
	\tfrac{\partial \varrho_1}{\partial s}\right)=
	\displaystyle\alpha
             \left(\tfrac{v}{r}\,  \varrho_2+
	\tfrac{B_2}{rB} w\,  \varrho_3\right)+F_1(\mathbf u,p),\vspace{2mm}\\
        \displaystyle  \varrho_2+\alpha\left(u \,
	 \tfrac{\partial \varrho_2}{\partial r}
	+\tfrac{v}{r} \, \tfrac{\partial \varrho_2}{\partial\theta}
	 +\tfrac{w}{B} \, 
	\tfrac{\partial \varrho_2}{\partial s}\right)=
	\displaystyle-\alpha
            \left(\tfrac{v}{r} \,  \varrho_1
	+\tfrac{B_1}{rB} w \, \varrho_3\right)
	+F_2(\mathbf u,p),\vspace{2mm}\\
	\displaystyle\varrho_3+\alpha\left(u \,
	 \tfrac{\partial \varrho_3}{\partial r}
	+\tfrac{v}{r} \, \tfrac{\partial \varrho_3}{\partial\theta}
	 +\tfrac{w}{B} \, 
	\tfrac{\partial \varrho_3}{\partial s}\right)
	=\displaystyle\alpha w
             \left(\tfrac{B_1}{rB} \,\varrho_2-
\tfrac{B_2}{rB}\,  \varrho_1\right)+F_3(\mathbf u,p),
           \end{array}\right.\end{equation}
where $\mathbf F_i(\mathbf u,p)=\left(\nabla\cdot\left(\mathbf L(\mathbf u)-\alpha
	p(\nabla\mathbf u)^T-{\cal R}e \, \mathbf u\otimes \mathbf u\right)
	\right)_i$.
We consider flows which are fully developed. The three components of the
velocity are then independent of the variable $s$, i.e.
	\begin{equation}\label{pres1}\tfrac{\partial u}{\partial s}=
	\tfrac{\partial v}{\partial s}=
	\tfrac{\partial w}{\partial s}\equiv0,\end{equation}
and consequently the axial component of the pressure gradient is a constant;
	\begin{equation}\label{prescurved} 
	\tfrac{\partial\pi}{\partial s}=-p^\ast.\end{equation}
Our aim now is to deduce an information on 
$\tfrac{\partial p}{\partial s}$.
 Taking into account the fact that 
	$$\pi=\left(I+\alpha\left(u\,\tfrac{\partial }{\partial r}
	+\tfrac{v}{r} \, \tfrac{\partial }{\partial \theta}+\tfrac{w}{B}\,
	\tfrac{\partial }{\partial s}\right)\right)p,$$
we deduce from (\ref{pres1})-(\ref{prescurved}) that
	$$\tfrac{\partial\pi}{\partial s}=-p^\ast=\left(I+\alpha\left(u\,\tfrac{\partial }{\partial r}
	+\tfrac{v}{r} \, \tfrac{\partial }{\partial \theta}+\tfrac{w}{B}\,
	\tfrac{\partial }{\partial s}\right)\right)\tfrac{\partial p}{\partial s},$$
and thus,
	$$\left(I+\alpha\left(u\,\tfrac{\partial }{\partial r}
	+\tfrac{v}{r} \, \tfrac{\partial }{\partial \theta}+\tfrac{w}{B}\,
	\tfrac{\partial }{\partial s}\right)\right)  \left(\tfrac{\partial p}{\partial s}
	+p^\ast\right)=0.$$
Since the operator $I+\alpha\left(u\,\tfrac{\partial }{\partial r}
	+\tfrac{v}{r} \, \tfrac{\partial }{\partial \theta}+\tfrac{w}{B}\,
	\tfrac{\partial }{\partial s}\right)$ is invertible, we obtain
	\begin{equation}\label{press2}\tfrac{\partial p}{\partial s}=-p^\ast.\end{equation}
A direct consequence of (\ref{pres1}), (\ref{press2}) and $(\ref{stokescurv})_{1,2,3}$ is that the auxiliary unknown $\boldsymbol\varrho$ is independent of $s$, and
	\begin{equation}\label{sig}\tfrac{\partial \varrho_1}{\partial s}=
	\tfrac{\partial \varrho_2}{\partial s}=
	\tfrac{\partial \varrho_3}{\partial s}\equiv 0.\end{equation}
Taking into account (\ref{pres1}), (\ref{press2})-(\ref{sig}) and replacing in system 
(\ref{stokescurv})-(\ref{transportcurv}), we easily see that the problem of fully developed 
second grade fluid defined in 
$$
\Sigma=\left\{(r,\theta)\in \mathbb R^2\mid \, 0<r<1,\, 0\leqslant \theta <2\pi)\right\},
$$
 reads as
\begin{equation}
\left\{
\begin{array}{ll} 
\displaystyle	{\cal A}u+\tfrac{2}{r^{2}}\tfrac{\partial v}{\partial \theta}
	+\tfrac{B^2+B_2^2}{(rB)^2}u
	-\tfrac{B_1(B+B_2)}{(rB)^2}v
	+\tfrac{\partial p}{\partial r}=\varrho_1,\vspace{2mm}\\
\displaystyle{\cal A}v-\tfrac{2}{r^2}\tfrac{\partial u}{\partial \theta}
	+\tfrac{B_1}{(rB)^2}u+\tfrac{B^2+B_1^2}{(rB)^2}v+
	\tfrac{1}{r}\tfrac{\partial p}{\partial \theta}=\varrho_2,\vspace{2mm}\\
	\displaystyle{\cal A}w+\left(\tfrac
	{\delta}{B}\right)^2w-\tfrac{p^\ast}{B}=\varrho_3,\vspace{2mm}\\
 \displaystyle{\tfrac{\partial u}{\partial r}+\tfrac{1}{r}
	\tfrac{\partial v}{\partial \theta}+
	\tfrac{B+B_2}{rB}u-\tfrac{B_1}{rB}v=0,}
	  \end{array}\right.\label{stokescurv_fd}
\end{equation}
where 
$${\cal A}\zeta=-\tfrac{\partial ^2\zeta}{\partial r^2}-\tfrac{1}{r^2}
	\tfrac{\partial ^2\zeta}{\partial \theta^2}
	-\tfrac{B+B_2}{rB} 
	\tfrac{\partial \zeta}{\partial r}
	+\tfrac{B_1}{r^2B}
	\tfrac{\partial \zeta}{\partial \theta},$$
and 
\begin{equation}\label{transportcurv_fd}\left\{\begin{array}{lll}
	\displaystyle\varrho_1+\alpha\left(u \,
	 \tfrac{\partial \varrho_1}{\partial r}
	+\tfrac{v}{r} \, \tfrac{\partial \varrho_1}{\partial\theta}\right)=
	\displaystyle\alpha
             \left(\tfrac{v}{r}\,  \varrho_2+
	\tfrac{B_2}{rB} w\,  \varrho_3\right)+F_1(\mathbf u,p),\vspace{2mm}\\
        \displaystyle  \varrho_2+\alpha\left(u \,
	 \tfrac{\partial \varrho_2}{\partial r}
	+\tfrac{v}{r} \, \tfrac{\partial \varrho_2}{\partial\theta}\right)=
	\displaystyle-\alpha
            \left(\tfrac{v}{r} \,  \varrho_1
	+\tfrac{B_1}{rB} w \, \varrho_3\right)\displaystyle
	+F_2(\mathbf u,p),\vspace{2mm}\\
	\displaystyle\varrho_3+\alpha\left(u \,
	 \tfrac{\partial \varrho_3}{\partial r}
	+\tfrac{v}{r} \, \tfrac{\partial \varrho_3}{\partial\theta}\right)
	=\displaystyle\alpha
             \left(\tfrac{B_1}{rB} w \, \varrho_2-
\tfrac{B_2}{rB}w \,  \varrho_1\right)\displaystyle+F_3(\mathbf u,p).
           \end{array}\right.\end{equation}
We shall reformulate problem (\ref{stokescurv_fd})-(\ref{transportcurv_fd}) in a suitable weak form. To this end, we first dot-multiply both sides of the equations $(\ref{stokescurv_fd})_{1,2,3}$ by $(rB)^2\zeta$, $\zeta\in H_0^1(\Sigma)$, and integrate by parts over $\Sigma$. Observing that 
	$$\begin{array}{lll}
\left({\cal A}u,(rB)^2\zeta\right)
&=\displaystyle\left(-\tfrac{\partial ^2 u}{\partial r^2}-\tfrac{1}{r^2}
	\tfrac{\partial ^2 u}{\partial \theta^2}
	-\tfrac{B+B_2}{rB} 
	\tfrac{\partial u}{\partial r}
	+\tfrac{B_1}{r^2B}
	\tfrac{\partial u}{\partial \theta},(rB)^2\zeta\right)
\vspace{2mm}\\
&=\displaystyle\left(\tfrac{\partial u}{\partial r},\tfrac{\partial }{\partial r}((rB)^2\zeta)\right)+
\left(\tfrac{\partial u}{\partial \theta},\tfrac{\partial}{\partial \theta}(B^2\zeta)\right) -
\left(rB\tfrac{\partial u}{\partial r},(B+B_2)\zeta\right)+\left(B\tfrac{\partial u}{\partial \theta},B_1\zeta\right)
\vspace{2mm}\\
&=\displaystyle\left(rB\tfrac{\partial u}{\partial r},rB\tfrac{\partial \zeta}{\partial r}+(B+B_2)\zeta\right)+
\left(B\tfrac{\partial u}{\partial \theta},B\tfrac{\partial \zeta}{\partial \theta}-B_1\,\zeta\right) \equiv a\left(u,\zeta\right),
\end{array}$$
and 
	$$\left(\tfrac{\partial p}{\partial r},(rB)^2\zeta \right)
	=-\left(rB\,p,rB \tfrac{\partial \zeta}{\partial r}
	+2(B+B_2)\zeta\right)\equiv b_1\left(p,\zeta\right),$$
	$$\left(\tfrac{1}{r}\tfrac{\partial p}{\partial \theta},
	(rB)^2\zeta \right)=-\left(rB\,p,B \tfrac{\partial \zeta}{\partial \theta}
	-2B_1\,\zeta\right)\equiv b_2\left(p,\zeta\right),$$
we deduce that problem (\ref{stokescurv_fd}) can be written as follows
	\begin{equation}\label{stokescurv_fd3} 
	\left\{\begin{array}{ll} 
       a_i((u,v),\zeta)+b_i\left(p,\zeta\right)=
	\left({\sigma}_{i},
	\zeta\right)&\quad  \mbox{for all} \ \zeta\in  H_{0}^{1}(\Sigma) \qquad (i=1,2), \vspace{4mm}\\ 
	a_3(w,\zeta)
	=\left(r^2B \, p^\ast+{\sigma}_{3},
	\zeta\right) &\quad  \mbox{for all} \ \zeta\in H_{0}^{1}(\Sigma),
	\vspace{4mm}\\ 
	\left(\tfrac{\partial}{\partial r}\left(  rB u\right) 
	+\tfrac{\partial}{\partial\theta}(B v),\psi\right)=0&\quad 
	 \mbox{for all} \ \psi\in L^2_0(\Sigma),\end{array}
	\right.\end{equation} 
where  ${\sigma}_{i}=(rB)^2\boldsymbol\varrho_{i}$ and
	$$\begin{array}{ll} a_1((u,v),\zeta)= a(u,\zeta)+
	\left(\left(B^2+B_2^2\right)
	u-B_1(B+B_2)  v+2B^2 \,
	\tfrac{\partial v}{\partial\theta},\zeta\right),\vspace{2mm}\\
	a_2((u,v),\zeta)= a(v,\zeta)+\left(\left(B^2+B_1^2\right)
	v+B_1  u-2B^2 \, \tfrac{\partial u}{\partial\theta},
	\zeta\right), \vspace{2mm}\\
	a_3(w,\zeta)=a(w,\zeta)+\left(
	(r\delta)^{2}w,\zeta\right).\end{array}$$
 To deal with the transport equation, let us first observe that
  $$(rB)^ 3 u \tfrac{\partial \varrho_i}{\partial r}=
  rBu\tfrac{\partial {\sigma}_{i}}{\partial r}-2(B+B_2)u{\sigma}_{i},\qquad
(rB)^ 3 \tfrac{v}{r} \tfrac{\partial \varrho_i}{\partial \theta}=
  Bv\tfrac{\partial {\sigma}_{i}}{\partial \theta}-2B_1v{\sigma}_{i},$$
and thus
  $$\left(\displaystyle\varrho_i+\alpha\left(u \,
	 \tfrac{\partial \varrho_i}{\partial r}
	+\tfrac{v}{r} \, \tfrac{\partial \varrho_i}{\partial\theta}\right),
(rB)^3\tau\right)$$
	$$=\left({\sigma}_{i},
rB\tau\right)+\alpha\left(rBu\tfrac{\partial {\sigma}_{i}}{\partial r}
+Bv\tfrac{\partial {\sigma}_{i}}{\partial \theta},\tau\right)-2\left(\left((B+B_2)u+B_1v\right){\sigma}_{i},\tau\right) \qquad \mbox{for all}
	\ \tau\in L^2(\Sigma).$$
Multiplying both sides of (\ref{transportcurv_fd}) by $(rB)^3\tau$, we obtain 
\begin{equation}\label{transportcurv_fd3}
	\left(rB{\sigma}_{i},\tau\right)
	+\alpha\left(rBu\tfrac{\partial {\sigma}_{i}}{\partial r}
+Bv\tfrac{\partial {\sigma}_{i}}{\partial \theta},\tau\right)=
\left(G_i\left(\mathbf u,p,\boldsymbol\sigma\right),\tau\right), \end{equation}
 where 
$$\begin{array}{ll}G_1\left(\mathbf u,p,\boldsymbol\sigma\right)=(rB)^3F_1(\mathbf u,p)+\alpha \left(2\left((B+B_2)u -B_1 v\right)
	 {\sigma}_1+B v\, {\sigma}_2+
	B_2 w \, {\sigma}_3\right),\vspace{2mm}\\
G_2\left(\mathbf u,p,\boldsymbol\sigma\right)=(rB)^3F_2(\mathbf u,p)+\alpha
            \left(2(\left(B+B_2)u -B_1 v\right){\sigma}_2-B v \, {\sigma}_1
	-B_1 w \, {\sigma}_3\right),\vspace{2mm}\\
	G_3\left(\mathbf u,p,\boldsymbol\sigma\right)=(rB)^3F_3(\mathbf u,p)+\alpha
	  \left(2\left((B+B_2)u -B_1 v\right){\sigma}_3+\left(B_1  {\sigma}_2-B_2 
	     {\sigma}_1\right)w\right).
\end{array}
$$
\section{Numerical approximation}
\subsection{Setting of the approximated problem}
Let $\{\mathcal{T}_h\}_{h>0}$ be a family of regular triangulations defined over the
rectangle $\Sigma$. We consider the following finite element spaces 
	$$V_{h}=\left\{v_{h}\in C(\overline\Sigma)
	\cap H_0^1(\Sigma)\mid
	{v_{h}}_{\mid K}\in \mathbb{P}_{2}(K)\quad \mbox{for all} \ 
	 K\in\mathcal{T}_{h}
	\right\}, $$
	$${M}_{h}=
	\left\{q_{h}\in C(\overline\Sigma)\cap L_{0}^{2}(\Sigma) \mid
	{q_{h}}_{\mid K}\in \mathbb{P}_{1}(K)\quad \mbox{for all} \ 
	 K\in\mathcal{T}_{h}\right\},$$
	$$T_{h}=\left\{{\tau}_{h}\in L^2(\Sigma)\mid {{\tau}_{h}}_{\mid K}\in
	\mathbb{P}_{1}\quad \mbox{for all} \  K\in \mathcal{T}_h \right\}.$$
For each $K\in {\cal T}_h$, with boundary $\partial K$, the inflow edge ${\partial K}^-$ of $K$ is defined by $\partial K^-(\mathbf w)=\displaystyle\left\{s\in \partial K\mid \mathbf w\cdot
           (n_r,n_\theta)(s)<0\right\}$ where $(n_r,n_\theta)$ is the outward unit normal vector
 to element $K$. For $s\in \partial K$ such that $\mathbf w\cdot
           (n_r,n_\theta)(s)\neq 0$, we define the left hand and right hand limits $\tau^-$ and $\tau^+$ as 
         $$\tau^-(s)=\lim_{\varepsilon\rightarrow 0^-}\tau(s+
        \varepsilon \mathbf w(s)),\qquad \tau^+(s)=
	\lim_{\varepsilon\rightarrow 0^+}\tau(s+
        \varepsilon \mathbf w(s)).$$
System (\ref{stokescurv_fd3})-(\ref{transportcurv_fd3}) is approximated 
by the following coupled system \vspace{1mm}
$$\hspace{-2cm}\text{Find } (\mathbf u_h\equiv (u_{h}, v_h,w_h), p_h,
	\boldsymbol{\sigma}_{h})
 \in (V_h)^3\times M_h\times (T_h)^3 \text{ solution of }$$
 \begin{equation}\label{stokescurv2} 
\left\{\begin{array}{l}a_i\left((u_h,v_h),\zeta_h\right)+b_i\left(p_h,\zeta_h\right)=
	\left({\sigma}_{h,i},
	\zeta_h\right)\qquad  \mbox{for all} \ 
 \zeta_h\in  V_h \qquad (i=1,2), \vspace{4mm}\\ 
	a_3(w_h,\zeta_h)
	=\left(r^2B \, p^\ast+{\sigma}_{h,3},\zeta_h\right)\qquad  \mbox{for all} \ 
 \zeta_h\in V_h,	\vspace{4mm}\\ 
	\left(\tfrac{\partial}{\partial r}\left(  rB u_h\right) 
	+\tfrac{\partial}
	{\partial\theta}(B v_h),\psi_h\right)=0\qquad  \mbox{for all} \ 
 \psi_h\in M_h,
\end{array}
	\right. \vspace{2mm}
\end{equation} 
\begin{equation}\label{transport_fd4}\left(rB \, {\sigma}_{h,i}, \tau_h\right)+ {\cal B}_h\left(\alpha{u}_h,\alpha v_h,
 {\sigma}_{h,i},
	 \tau_h\right)=\left(G_i(\mathbf u_h,p_h,\boldsymbol{\sigma}_h),\tau_h\right) \qquad 
	\mbox{for all} \ \tau_h\in T_h,\end{equation}
where  ${\cal B}_h$ defined by
$$\begin{array}{ll}
	{\cal B}_h(u,v,\sigma, \tau)=&\displaystyle \left(rB 
	u \tfrac{\partial \sigma}{\partial r}+
	B v \tfrac{\partial \sigma}{\partial\theta},\tau\right)_h
	+\tfrac{1}{2}\left(\left(\tfrac{\partial }{\partial r}\left(rB u\right)
	+\tfrac{\partial }{\partial \theta}\left(Bv\right)\right) \sigma,\tau\right)\vspace{2mm}\\
	&\displaystyle-\langle \sigma^+
	-\sigma^-,\tau^+ \rangle_{h,(u,v)},\end{array}$$
with 
	$$\left(\zeta,\psi\right)_h=\displaystyle \sum_{K\in {\cal T}_h}
	 \int_K\zeta\psi\,drd\theta, $$
	$$\langle \sigma,\tau \rangle_{h,(u,v)}=
	\sum_{K\in\mathcal{T}_h}
	 \int_{\partial K^-(rBu,Bv)} 
	\sigma\,\tau\, 
	 (rBu\, n_r+Bv\, n_\theta)\, ds.$$
 By a standard integration by parts, we obtain
$$\begin{array}{ll}\displaystyle
	\alpha\left( rB u \tfrac{\partial \sigma}{\partial r}+
	B v \tfrac{\partial \sigma}{\partial\theta},\tau\right)_h=&
	\displaystyle-\alpha	\left( rB u \tfrac{\partial \tau}{\partial r}+
	B v \tfrac{\partial \tau}{\partial\theta},\sigma\right)_h
	-\alpha\left(\left(\tfrac{\partial }{\partial r}\left(rB u\right)
	+\tfrac{\partial }{\partial \theta}\left(Bv\right)\right)\sigma,\tau\right)\vspace{2mm}\\ 
	&\displaystyle+\langle \sigma^+,\tau^+ \rangle_{h,(\alpha u,\alpha v)}-\langle \sigma^-
	,\tau^- \rangle_{h,(\alpha u,\alpha v)},\end{array}$$
and deduce that the bilinear form ${\cal B}_h$ satisfies
$$\begin{array}{ll}
	{\cal B}_h(\alpha u,\alpha v,\sigma, \tau)&=-\alpha\left( rB 
	u \tfrac{\partial \tau}{\partial r}+
	B v \tfrac{\partial \tau}{\partial\theta},\sigma\right)_h
	-\tfrac{\alpha}{2}\left(\left(\tfrac{\partial }{\partial r}\left(rB u\right)
	+\tfrac{\partial }{\partial \theta}\left(Bv\right)\right)\tau,\sigma\right)\vspace{2mm}\\
	&+\langle \sigma^-,\tau^+-\tau^- \rangle_{h,(\alpha u,\alpha v)}.\end{array}$$
Therefore,  problem (\ref{transport_fd4}) can be equivalently written as    
\begin{align}\label{transport_fd5}&\left(
 {\sigma}_{h,i}, rB\tau_h-\alpha \left(rB 
	u \tfrac{\partial \tau_h}{\partial r}+
	B v \tfrac{\partial \tau_h}{\partial\theta}\right)\right)
-\tfrac{\alpha}{2}\left(
 {\sigma}_{h,i}, \left(\tfrac{\partial }{\partial r}\left(rB u\right)
	+\tfrac{\partial }{\partial \theta}\left(Bv\right)\right)\tau_h\right)\nonumber\\
	&+\displaystyle\langle
	{\sigma}_{h,i}^-,\tau_h^+-\tau_h^-\rangle_{h,(\alpha u,\alpha v)}=\left(G_i(\mathbf u_h,p_h,\boldsymbol{\sigma}_h),\tau_h\right)\qquad \mbox{for all} \ \tau_h\in T_h.\end{align}
\subsection{Algorithm}
The numerical algorithm to solve the approximated problem 
(\ref{stokescurv2})-(\ref{transport_fd5}) is based on Newton's method, with the non-linear part 
explicitely calculated at each iteration step. As indicated below, at each step of the linearization process, a Stokes system is
solved for $(\mathbf{u},p)$, a Poisson equation is solved for the axial 
velocity $w$ and a transport equation for 
 $\boldsymbol\sigma$. \vspace{2mm}\\
$\bullet$ Given an iterate $\boldsymbol{\sigma}_h^k\equiv
 ({\sigma}_{h,1}^k,{\sigma}_{h,2}^k,{\sigma}_{h,3}^k)$, find 
$\mathbf{u}_h^k\equiv \left(u^k_h,v^k_h,w^k_h\right)$ and
$p^k_h$ solutions of the Stokes system 
$$ ({\cal P}_k)\qquad	\left\{\begin{array}{l}
       a_i\left((u^k_h,v^k_h),\zeta_h\right)+b_i\left(p_h^k,\zeta_h\right)=
	\left({\sigma}_{h,i}^k,
	\zeta_h\right)\qquad  \mbox{for all} \ 
 \zeta_h\in  V_h \qquad (i=1,2)\vspace{3mm}\\ 
	a_3\left(w_h^k,\zeta_h\right)
	=\left(r^2B \, p^\ast+{\sigma}_{h,3}^k,\zeta_h\right)\qquad  \mbox{for all} \ 
 \zeta_h\in V_h
	\vspace{3mm}\\ 
	\left(\tfrac{\partial}{\partial r}\left(  rB u_h^k\right) 
	+\tfrac{\partial}
	{\partial\theta}(B v_h^k),\psi_h\right)=0\qquad  \mbox{for all} \ 
 \psi_h\in M_h 
	\end{array}
	\right.
$$
$\bullet$ Calculate the new iterate $\boldsymbol\sigma^{k+1}_h$ as the solution of the following transport problem 
$$\begin{array}{ll}
\left(rB \tau_{h}-\alpha \left( rB 
	u^k_h \tfrac{\partial \tau_h}{\partial r}+
	B v^k_h \tfrac{\partial \tau_h}{\partial\theta}\right) ,  
  {\sigma}_{h,i}^{k+1}\right)
-\tfrac{\alpha}{2}\left( \left(\tfrac{\partial }{\partial r}\left(rB u_h^k\right)
	+\tfrac{\partial }{\partial \theta}\left(Bv_h^k\right)\right)\tau_h,  
 {\sigma}_{h,i}^{k+1}\right)\vspace{2mm}\\
\displaystyle +\left\langle
	\left({\sigma}_{h,i}^{k+1}\right)^-,\tau_h^+-\tau_h^-\right\rangle_{h,(\alpha u^k_h,\alpha v^k_h)}=\left(G_i(\mathbf u^k_h,p^k_h,
\sigma^{k}_h)
,\tau_h\right) \qquad \mbox{for all}  \ \tau_h\in T_h
\end{array}$$
$\bullet \ $ Find $(\mathbf{u}_h^{k+1},p_h^{k+1})$ solution
 of the Stokes system $({\cal P}_{k+1})$.

\section{Numerical results}
\label{S7} 

\noindent In order to study the non-Newtonian effects of the flows, we compare
the quantitative and qualitative behaviour of the axial velocity and  
the stream function of both creeping and inertial flows of second-grade fluids.
 Using a continuation method on the characteristic parameters 
(the Reynolds number ${\cal R}e$ and the viscoelastic parameter $\alpha$),
 we obtain numerical results in 
different flow situations. Introducing an adimensionalised stream function $\psi$ we obtain the
following identities for the velocity, \vspace{-2mm}
$$
u=-\tfrac{1}{rB}\tfrac{\partial\psi}{\partial\theta}\ ,
\qquad
v=\tfrac{1}{B}\tfrac{\partial\psi}{\partial r}\ .
$$
%
\subsection{Newtonian flows}  
\begin{figure}[!ht]
\centering
\begin{tabular}{cc}
\includegraphics[width=48mm]{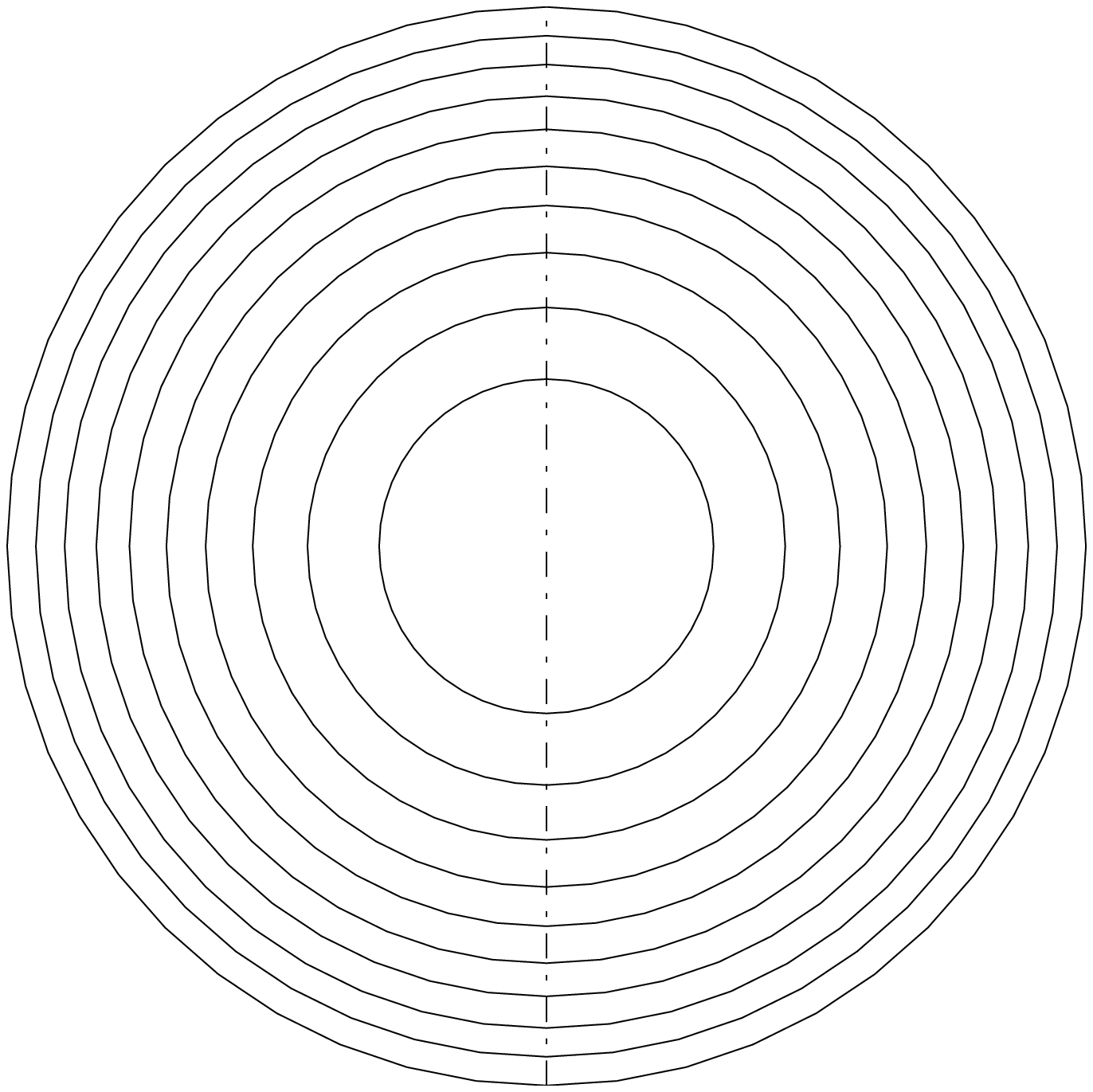} & \hspace{7mm}
\includegraphics[width=48mm]{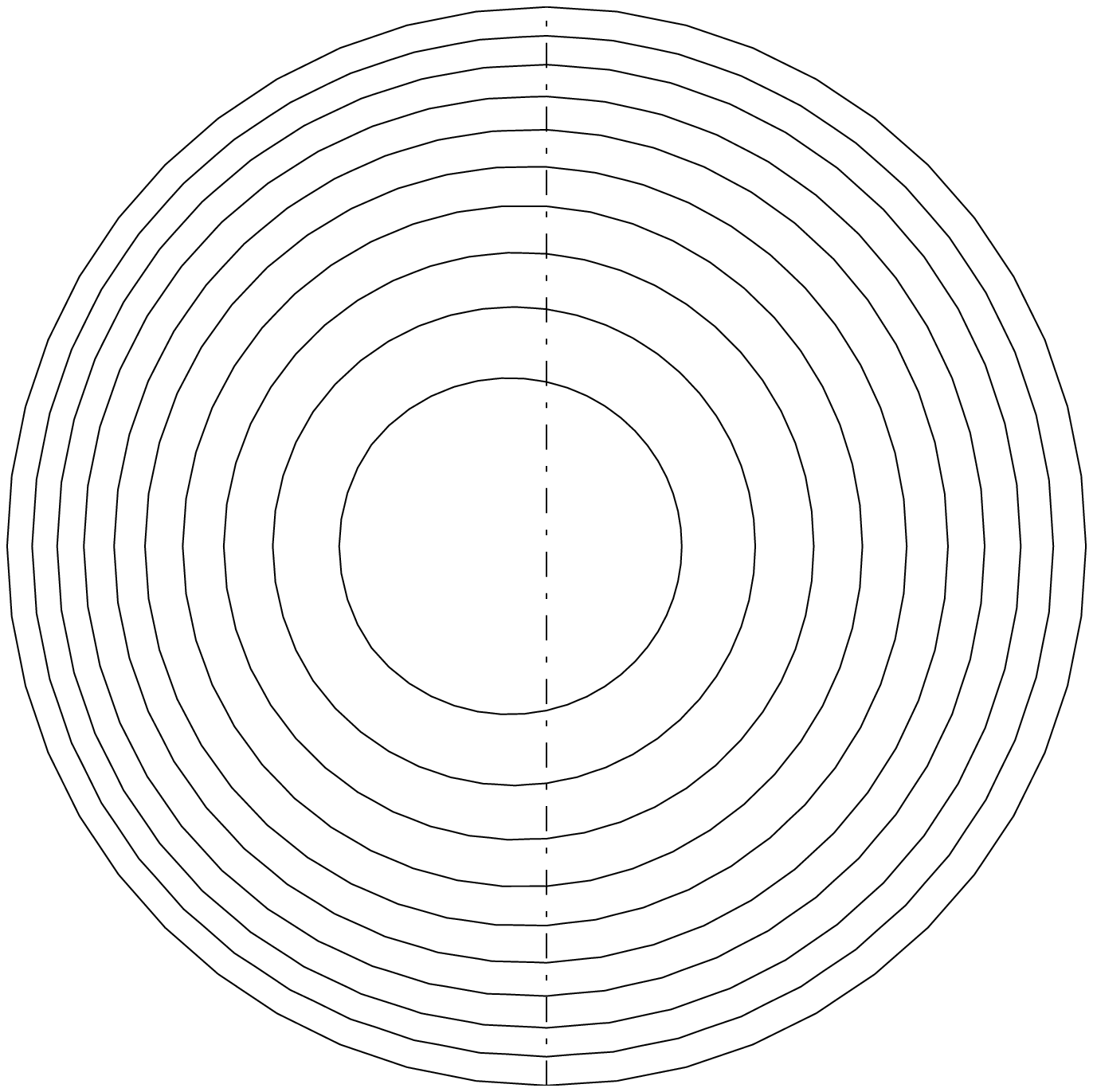} \\
{\small (a)$\ \delta=0.001$} & \hspace{10mm} {\small (b) $\ \delta=0.2$} 
\end{tabular}
\caption{\small \mbox{Qualitative behaviour of the axial velocity for creeping Newtonian flow}.}
\label{poiseuille}
\end{figure}

%
It is interesting to compare the qualitative behaviour of the flow for 
second-grade fluids with that of Newtonian fluids.
 For this purpose, we first consider typical contours of the axial velocity
and of the stream function in the case of Newtonian fluids.
In the case of creeping fluids (${\cal R}e\!=\!0$),
 a Poiseuille solution 
is displayed for a small curvature ratio ($\delta\!=\!0.001$).
There is no secondary motion and the contours of the 
axial velocity $w$ are circles, centered about the central axis (see 
Figure $\ref{poiseuille}_a$). In contrast, as can be seen in Figure 
$\ref{poiseuille}_b$, the contours are shifted away from the center towards the 
inner wall when the curvature ratio increases.

\begin{figure}[!ht]
\SetLabels
\L (0.48*0.986) $\leftarrow$ \\
\L (0.48*0.02) $\leftarrow$ \\
\endSetLabels
\leavevmode
\centering
\begin{tabular}{cc}
\strut\AffixLabels{
\includegraphics[width=48mm]{./figure03_1.eps}} & \hspace{7mm}
\includegraphics[width=48mm]{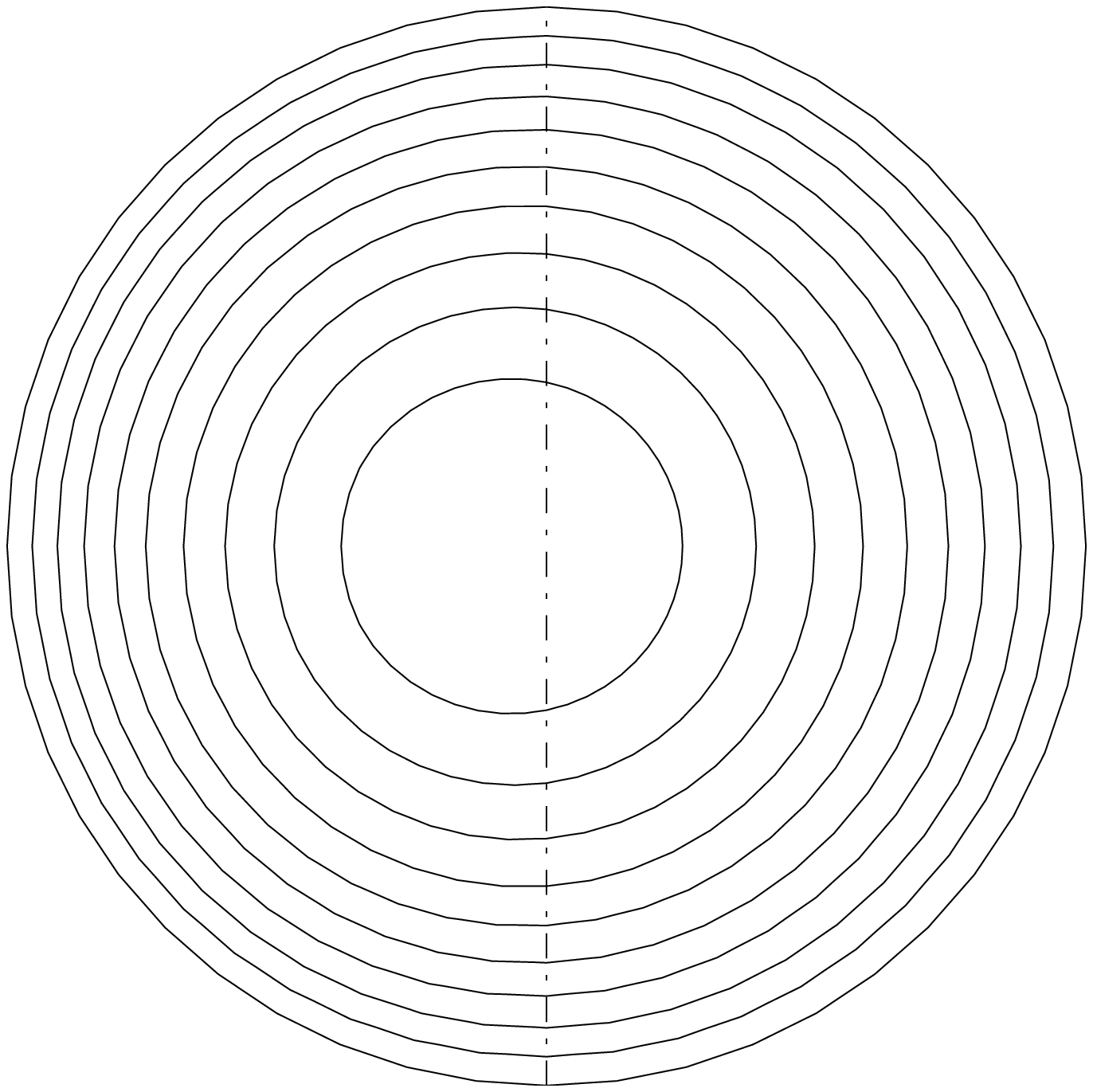}
\end{tabular}
\caption{\small Qualitative behaviour of the streamlines and the axial velocity for Newtonian flows with
 $\mathcal{R}e\!=\!5$ and $\delta\!=\!0.2$.}
\label{NS-fig2}
\end{figure}
\noindent For inertial Newtonian fluids 
(${\cal R}e\!\neq\! 0$) a slight curvature of the pipe axis 
induces centrifugal forces on the fluid 
and consequently secondary flows, sending fluid outward along the symmetry axis 
and returning along the upper and lower curved surfaces. A pair of symmetric 
vortices is then superposed to the Poiseuille flow. 
This can be seen in Figure $\ref{NS-fig2}$, where we contours of
 the axial velocity and the stream-function are presented for 
${\cal R}e=5$ and $\delta=0.2$. \vspace{2mm}\\
\noindent The solutions obtained with FEM and 
Robertson's perturbation method show a good agreement, 
in accordance with the predictions. In particular, for this value of 
the curvature ratio, we can observe a shift from the center. 
The secondary flow is
counter-clockwise in the upper half of the cross--section and clockwise in the
lower half.
\begin{figure}[!ht]
\SetLabels
\L (0.48*0.830) {\small{$\psi_{\max}$}} \\
\L (0.48*0.200) {\small{$\psi_{\min}$}} \\
\L (0.98*0.00) {\tiny{$\mathcal{R}\! e$}} \\
\endSetLabels
\leavevmode
\centering
\begin{tabular}{c}
$$
\strut\AffixLabels{
\includegraphics[width=45mm,angle=270]{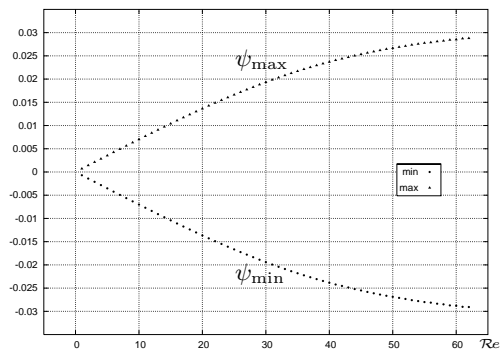}}
$$  \vspace{1mm}
\end{tabular}
\caption{\small Profile of the extrema values of the stream function for inertial Newtonian flows.}
\label{extrema_Newt}
\end{figure}
\noindent As already known, the streamlines corresponding to the Newtonian flows show a qualitative 
and quantitative symmetry relatively to the axis $r=0$.  To point out this phenomenon, and for 
comparison purposes, we first plot the profile of the extrema values of the stream function for 
different Reynolds numbers. As can be seen in Figure \ref{extrema_Newt}, the modulus of the 
extrema increase with Reynolds numbers. 

%
\subsection{Creeping viscoelastic flows}  
%
In this section, we are interested in the viscoelastic behaviour of 
the fluid in the case of creeping flows, and especially in 
the behaviour of the secondary motions. Our numerical results (using the FEM) indicate 
changes in the flow characteristics and suggest that at zero Reynolds number, 
fluid viscoelasticity promotes a secondary flow.
\begin{figure}[!ht]
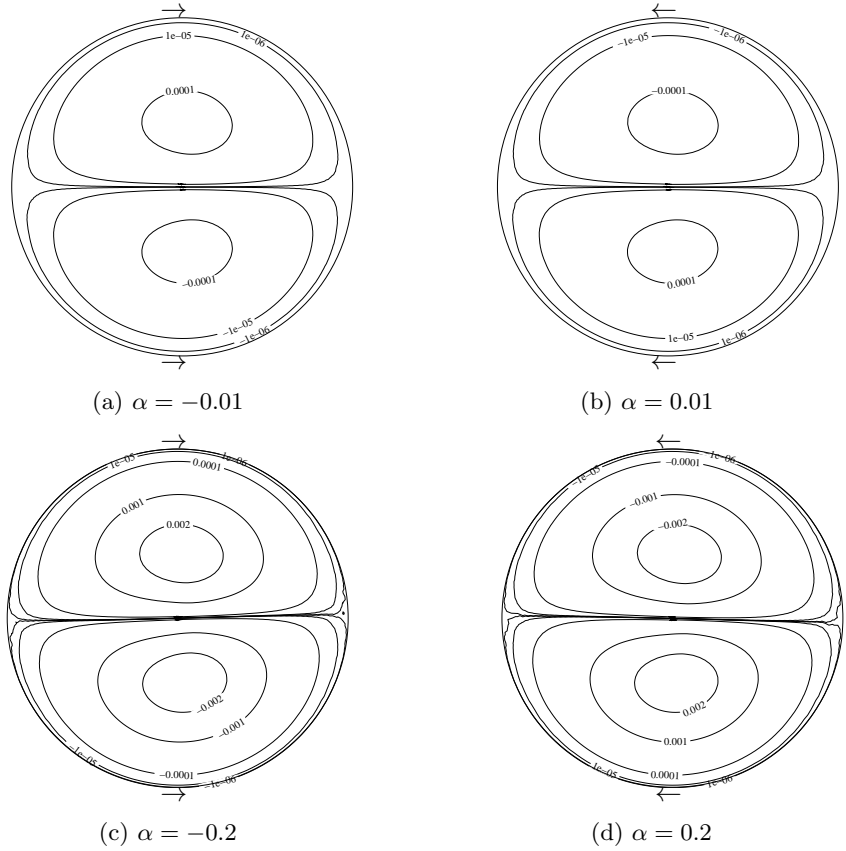

\SetLabels
\L (0.48*0.986) $\rightarrow$ \\
\L (0.48*0.02) $\rightarrow$ \\
\endSetLabels
\leavevmode
\centering
\begin{tabular}{cc}
$$
\strut\AffixLabels{
\includegraphics[width=48mm]{./figure05_1.eps}}
$$
\SetLabels
\L (0.48*0.986) $\leftarrow$ \\
\L (0.48*0.02) $\leftarrow$ \\
\endSetLabels
\leavevmode
& \hspace{7mm}
$$
\strut\AffixLabels{
\includegraphics[width=48mm]{./figure05_2.eps} 
}
$$ \\
{\small (a)$\ \alpha=-0.01$} & \hspace{4mm} {\small (b)$\ \alpha=0.01$}\vspace{3mm}
\end{tabular}

\SetLabels
\L (0.48*0.986) $\rightarrow$ \\
\L (0.48*0.02) $\rightarrow$ \\
\endSetLabels
\leavevmode
\centering
\begin{tabular}{cc}
$$
\strut\AffixLabels{
\includegraphics[width=48mm]{./figure05_3.eps} 
}
$$
\SetLabels
\L (0.48*0.986) $\leftarrow$ \\
\L (0.48*0.02) $\leftarrow$ \\
\endSetLabels
\leavevmode
& \hspace{7mm}
$$
\strut\AffixLabels{
\includegraphics[width=48mm]{./figure05_4.eps} 
}
$$
\\
{\small (c)$\ \alpha=-0.2$} & \hspace{4mm} {\small (d)$\ \alpha=0.2$}
\end{tabular}
\caption{\small Qualitative behaviour of the streamlines and the axial velocity for creeping viscoleastic 
flows with $\delta\!=\!0.2$.
}
\label{creeping3}
\end{figure}

\vspace{2mm}

\begin{figure}[!ht]
\SetLabels
\L (0.48*0.850) {\small{$\psi_{\max}$}} \\
\L (0.48*0.150) {\small{$\psi_{\min}$}} \\
\L (0.98*-0.01) {\tiny{$\alpha$}} \\
\endSetLabels
\leavevmode
\centering
\begin{tabular}{cc}
$$
\strut\AffixLabels{
\includegraphics[width=38mm,angle=270]{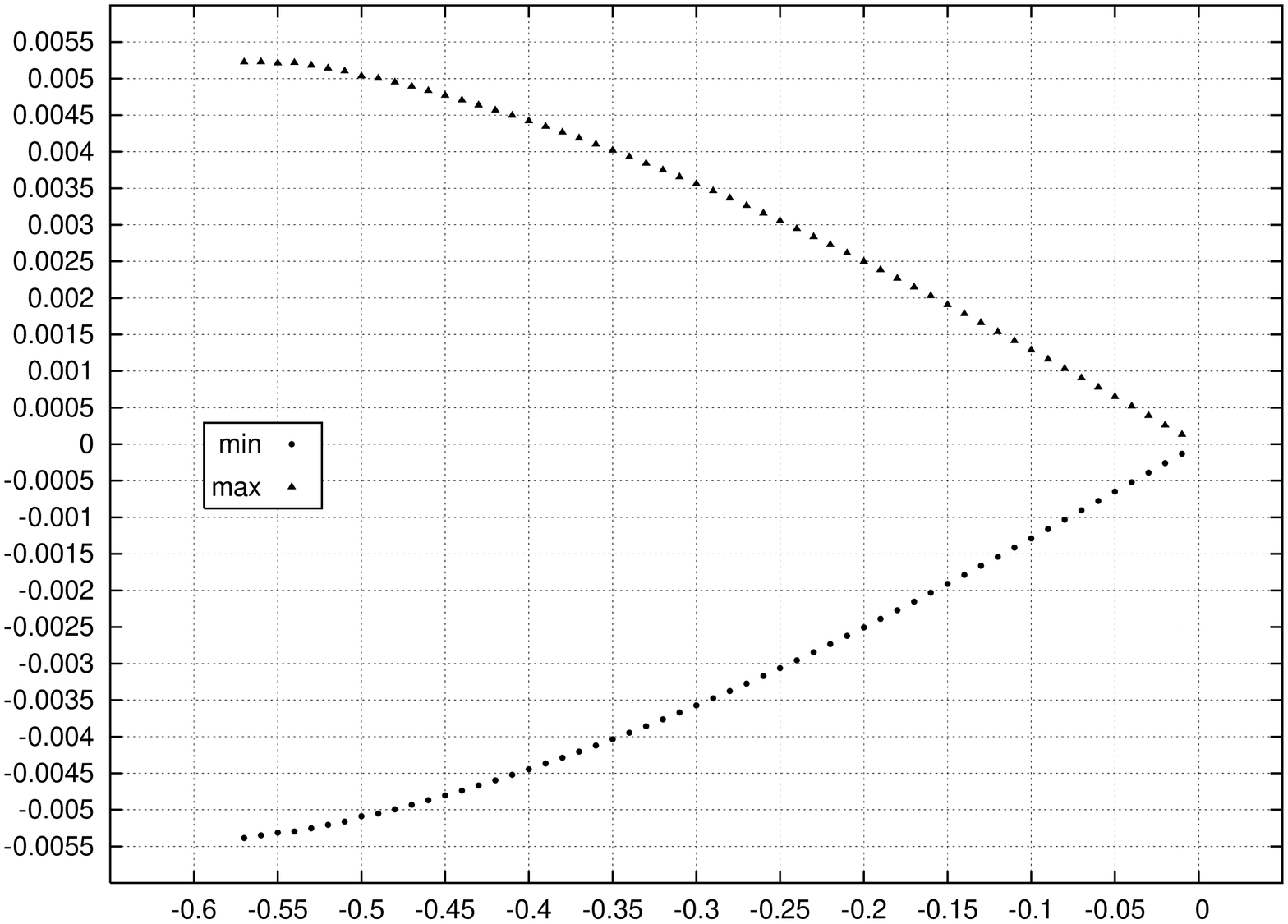}
}
$$ \
\SetLabels
\L (0.48*0.850) {\small{$\psi_{\max}$}} \\
\L (0.48*0.150) {\small{$\psi_{\min}$}} \\
\L (0.98*-0.01) {\tiny{$\alpha$}} \\
\endSetLabels
\leavevmode 
& \hspace{7mm} 
$$
\strut\AffixLabels{
\includegraphics[width=38mm,angle=270]{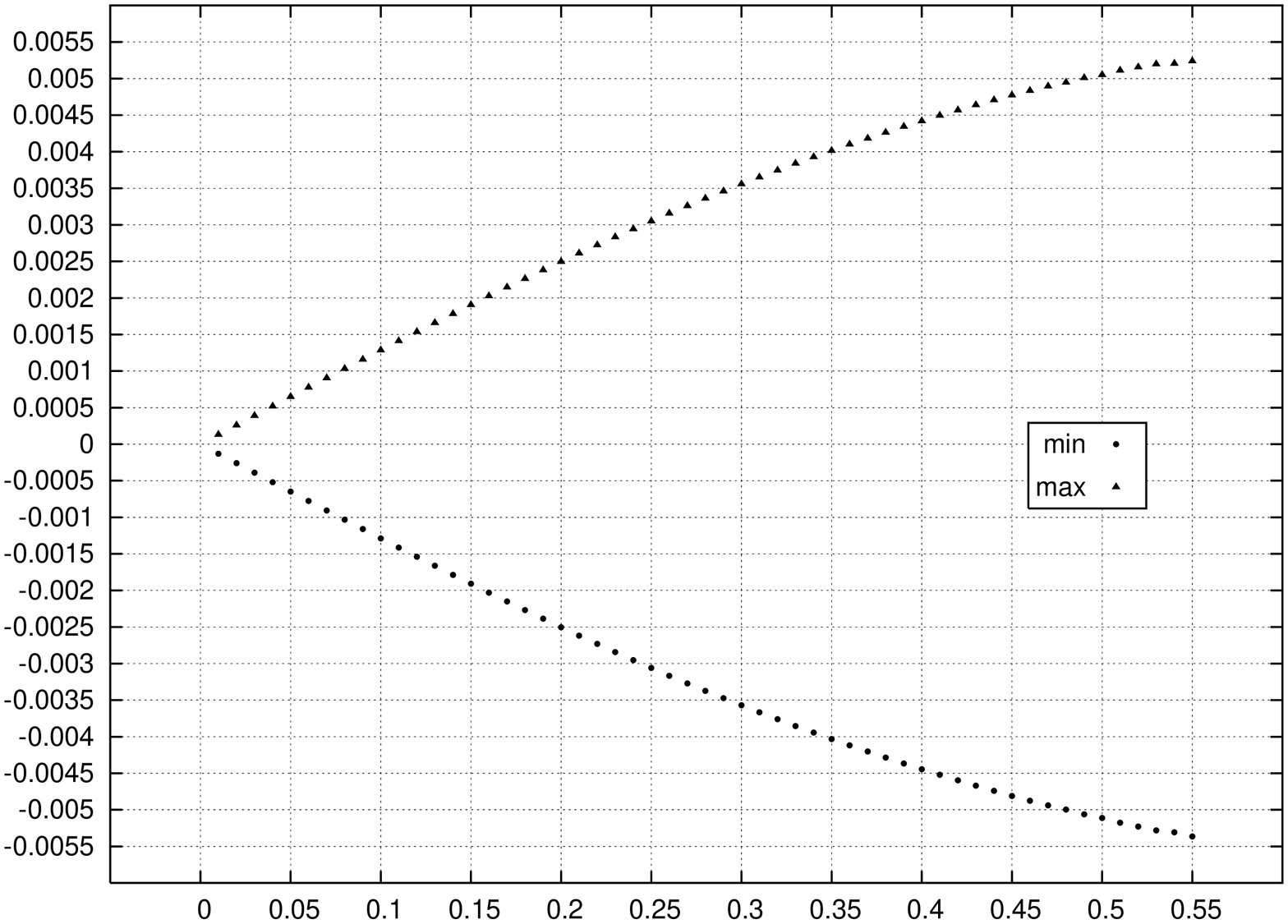} 
}
$$
 \vspace{1mm}\\
{\small (a)$\ \alpha$ negative} & 
{\small (b) $\ \alpha$ positive}
\end{tabular}
\caption{\small Profile of the extrema values of the stream function for creeping viscoelastic flows.}
\label{psi_extrema_creep}
\end{figure}

\noindent In a first step, we consider creeping fluids that are compatible with thermodynamics in the sense that 
the motions meet the Clausius-Duhem inequality (the viscoelastic parameter $\alpha$ is nonnegative).
In Figure  $\ref{creeping3}_b$ and Figure 
$\ref{creeping3}_d$, 
the streamlines are plotted for $\alpha=0.01$ and $\alpha=0.2$. 
As in the Newtonian case, the secondary flows are characterized by two 
counter-rotating vortices and their magnitude increase with
the elastic level (see Figure $\ref{psi_extrema_creep}_b$). The only difference is that
the flows remain symmetric relatively to an axis which is identical to the pipe 
centerline for $\alpha$ small and slightly rotates clockwise when
$\alpha$ increases (see Figures $\ref{creeping3}_b$ and $\ref{creeping3}_d$).

\noindent In a second step, we consider creeping fluids that are not compatible 
with the thermodynamics ($\alpha$ is negative). As in the positive case, the 
secondary flows involve non-zero values whose magnitude increase with $|\alpha|$ 
(see Figure $\ref{psi_extrema_creep}_a$). However, their behaviour is not similar. 
In particular, the orientation of the contours of the stream function is opposite. 
Moreover, the flow is symmetric relatively to an axis which is identical to the 
pipe centerline for $|\alpha|$ small, and slightly rotates counter-clockwise when 
$|\alpha|$ increases 
(see Figures $\ref{creeping3}_a$ and $\ref{creeping3}_c$).

\subsection{Inertial viscoelastic flows}  \vspace{-4pt}
%
Now, we are concerned with second-grade fluids
where the Reynolds number is non-zero, and with the combined effect of inertia and viscoelasticity.

\begin{figure}[!ht]
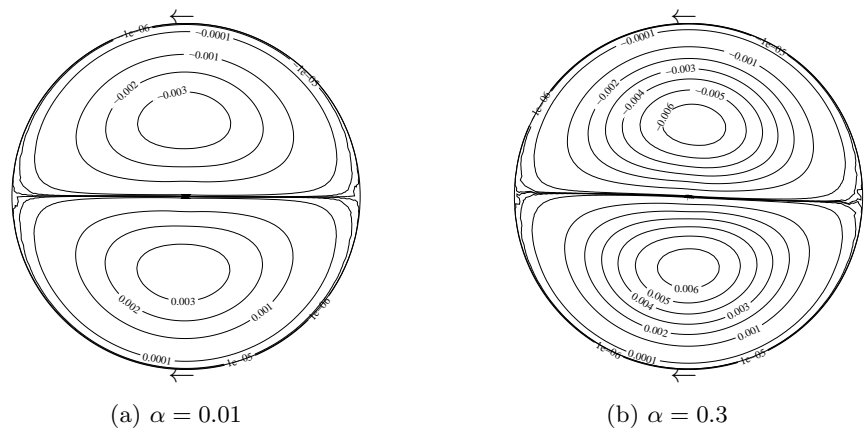

\SetLabels
\L (0.48*0.986) $\leftarrow$ \\
\L (0.48*0.02) $\leftarrow$ \\
\endSetLabels
\leavevmode
\centering
\begin{tabular}{cc}
$$
\strut\AffixLabels{
\includegraphics[width=49mm]{./figure07_1.eps} 
}
$$
\SetLabels
\L (0.48*0.986) $\leftarrow$ \\
\L (0.48*0.02) $\leftarrow$ \\
\endSetLabels
\leavevmode 
& \hspace{7mm}
$$
\strut\AffixLabels{
\includegraphics[width=49mm]{./figure07_2.eps} 
}
$$ \\
{\small (a)$\ \alpha=0.01$} & \hspace{4mm} {\small (b)$\ \alpha=0.3$}
\end{tabular}
\caption{\small Effect of the viscoelasticity on the streamlines with ${\mathcal R}e=5$ and  $\delta=0.2$.}
\label{visco-curve-pos1}
\end{figure}  
\begin{figure}[!ht]
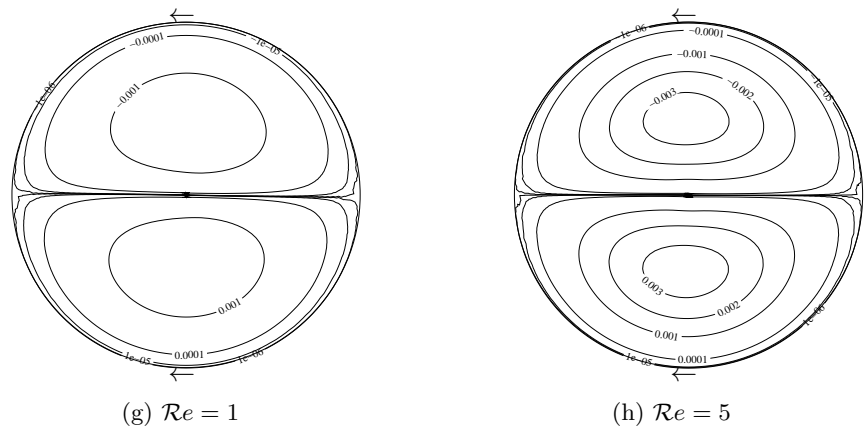

\SetLabels
\L (0.48*0.986) $\leftarrow$ \\
\L (0.48*0.02) $\leftarrow$ \\
\endSetLabels
\leavevmode
\centering
\begin{tabular}{cc}
$$
\strut\AffixLabels{
\includegraphics[width=49mm]{./figure08_1.eps} 
}
$$
\SetLabels
\L (0.48*0.986) $\leftarrow$ \\
\L (0.48*0.02) $\leftarrow$ \\
\endSetLabels
\leavevmode
& \hspace{7mm}
$$
\strut\AffixLabels{
\includegraphics[width=49mm]{./figure08_2.eps} 
}
$$ \\
{\small{ (g) $\mathcal{R}e=1$}} & \hspace{4mm} {\small{ (h) $\mathcal{R}e=5$}}
\end{tabular}
\caption{\small Effect of the inertia on the streamlines with $\alpha=0.1$ and  $\delta=0.2$.
}
\label{inert-curve-pos05}
\end{figure} \vspace{2mm}

\noindent As in the previous section, we first consider the case where the 
thermodynamics are satisfied. Setting the Reynolds to $5$ and increasing the 
viscoelastic parameter, we observe that the inertial and creeping viscoelastic 
flows have the same
behaviour: {\it globally Newtonian}, with lost of the symmetry relatively to the 
horizontal axis (see Figure \ref{visco-curve-pos1}). This symmetry is recovered 
if we fix the viscoelastic parameter and increase the Reynolds number 
(see Figure \ref{inert-curve-pos05}). 
Also, the magnitude of the secondary flows increase with $\alpha$ (see Figure
$\ref{extrema_viscoinert}_b$). \vspace{2mm}

\begin{figure}[!ht]
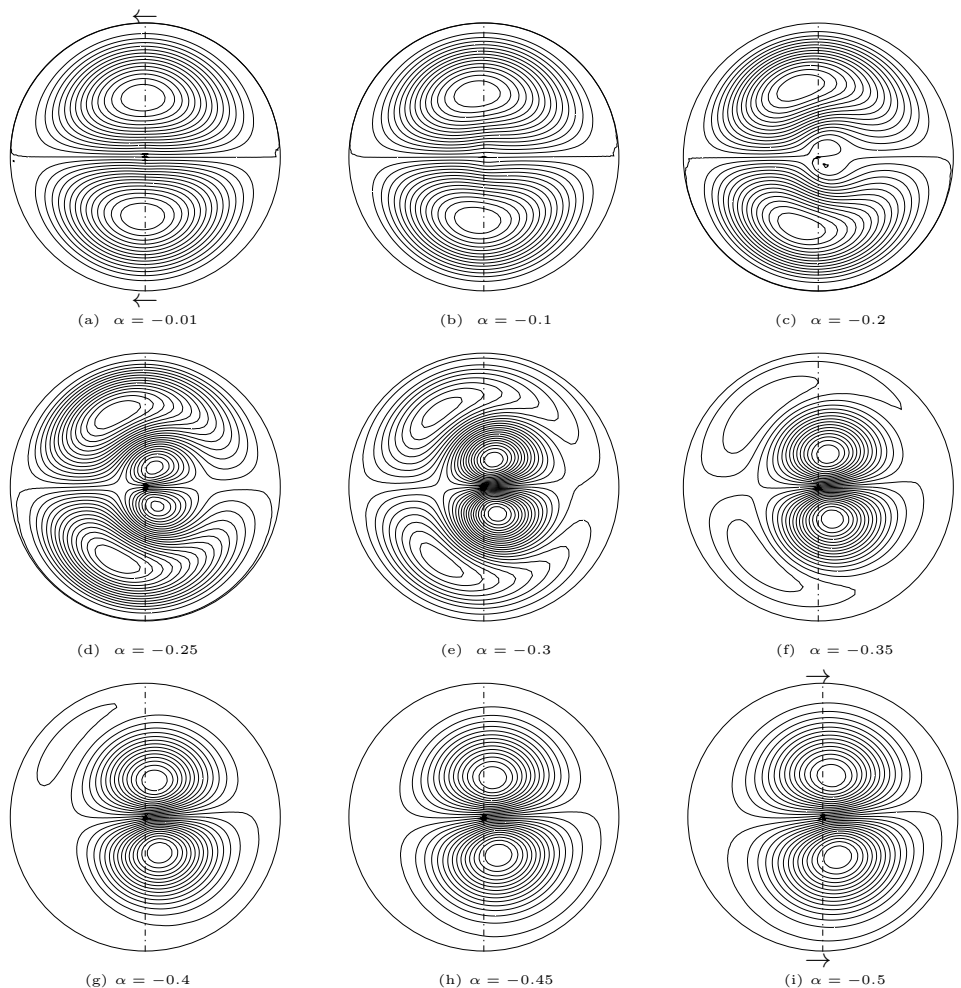

\SetLabels
\L (0.48*0.986) $\leftarrow$ \\
\L (0.48*0.0) $\leftarrow$ \\
\endSetLabels
\leavevmode
\centering
\begin{tabular}{ccc}
$$
\strut\AffixLabels{
\includegraphics[width=38mm]{./figure09_1.eps} 
}
$$
\SetLabels
\L (0.48*0.986) $\rightarrow$ \\
\L (0.48*0.0) $\rightarrow$ \\
\endSetLabels
\leavevmode
& \includegraphics[width=38mm]{./figure09_2.eps} 
& \includegraphics[width=38mm]{./figure09_3.eps}
\vspace{-2mm}\\
{\tiny (a) $\ \alpha=-0.01$} & \hspace{4mm} {\tiny (b) $\ \alpha=-0.1$}
& \hspace{4mm} {\tiny (c) $\ \alpha=-0.2$} \vspace{2mm}\\
\includegraphics[width=38mm]{./figure09_4.eps} 
& \includegraphics[width=38mm]{./figure09_5.eps} 
& \includegraphics[width=38mm]{./figure09_6.eps}
\vspace{-2mm}\\
{\tiny (d) $\ \alpha=-0.25$ }& \hspace{4mm} {\tiny (e) $\ \alpha=-0.3$} 
& \hspace{4mm} {\tiny (f) $\ \alpha=-0.35$}\vspace{2mm}\\
\includegraphics[width=38mm]{./figure09_7.eps} &
\includegraphics[width=38mm]{./figure09_8.eps} & 
$$
\strut\AffixLabels{
\includegraphics[width=38mm]{./figure09_9.eps}}
$$ \vspace{-2mm}\\
{\tiny (g)$\ \alpha=-0.4$} & \hspace{4mm} {\tiny (h)$\ \alpha=-0.45$}
& \hspace{4mm} {\tiny (i)$\ \alpha=-0.5$}
\end{tabular}
\caption{\small Effects of the viscoelasticity on the streamlines with ${\mathcal R}e=5$ and  $\delta=0.2$ with $\alpha$ negative.}
\vspace{3mm}
\label{visco-curve1}
\end{figure}  

\noindent The case of negative viscoelastic parameters is more complex.
When $\alpha$ is negative and small ($|\alpha|=0.01$), the nature and magnitude 
of the flow is qualitatively identical to that of a Newtonian fluid, 
the inertial effect being more pronounced and no effect due
 to viscoelasticity appears. When $|\alpha|$ increases, the solution obtained by 
 FEM shows reversal flows whose nature is very close to that obtained in the case 
 considered in the previous section when the secondary motion was generated by 
 fluid viscoelasticity.\vspace{1mm}

\noindent Figure \ref{visco-curve1} $\;$presents the streamlines for different 
values of $\alpha$ in the interval $[-0.01;-0.5]$.
For  $\alpha=-0.01$ we are still in the Newtonian regime. Two counter-rotating 
vortices appear, and the streamlines in the core region are more dense than elsewhere.
 However, in Figure $\ref{visco-curve1}_b$, a slight modification 
is already visible with the displacement of the vortices towards the inner wall.
 As $|\alpha|$ increases, these facts become more pronounced. 
 The values between $\alpha=-0.3$ and $\alpha=-0.2$ seems to be critical. The distortion of the 
streamlines is even more dramatic, and in addition to the vortices
 described below, another couple of vortices appear in the core region (see 
Figure $\ref{visco-curve1}_c$). 
Interestingly, the values in this part of the cross-section have 
an opposite sign to the ones near the boundaries, suggesting that the core is
 the transition region and that the reversal from one state to another 
initiates there.
This is confirmed by the results obtained for the cases where $\alpha=-0.3, -0.35, -0.4$,
showing that the reversal secondary flows grow around the new couple of
 vortices
 and that the changes occur from the core region 
to the regions near the boundary. It is clear that the size and strength of 
the new pair of vortices is more important, while a  slackening of the 
streamlines and vortices corresponding to the Newtonian state
 is observed.
In Figure $\ref{visco-curve1}_g$ the transition from one regime to the other
 is completed.
Finally, for $\alpha=-0.45$ due to the combined effect of curvature
and viscoelasticity, we observe a slight displacement of the center of the 
vortices to the outer wall. \noindent The same modifications can be observed in the 
quantitative bahaviour of the stream function. In particular, as can be seen in 
Figure $\ref{extrema_viscoinert}_a$, the modulus of the corresponding extremum values 
decrease for $|\alpha|\in [0,0.25]$, achieves a critical value for 
$|\alpha|\in [0.25,0.3]$ and increase for $|\alpha|\geq 0.3$.  

\begin{figure}[!ht]
\SetLabels
\L (0.38*0.780) {\small{$\psi_{\max}$}} \\
\L (0.38*0.200) {\small{$\psi_{\min}$}} \\
\L (0.98*-0.01) {\tiny{$\alpha$}} \\
\endSetLabels
\leavevmode
\centering
\begin{tabular}{cc}
$$
\strut\AffixLabels{
\includegraphics[width=38mm,angle=270]{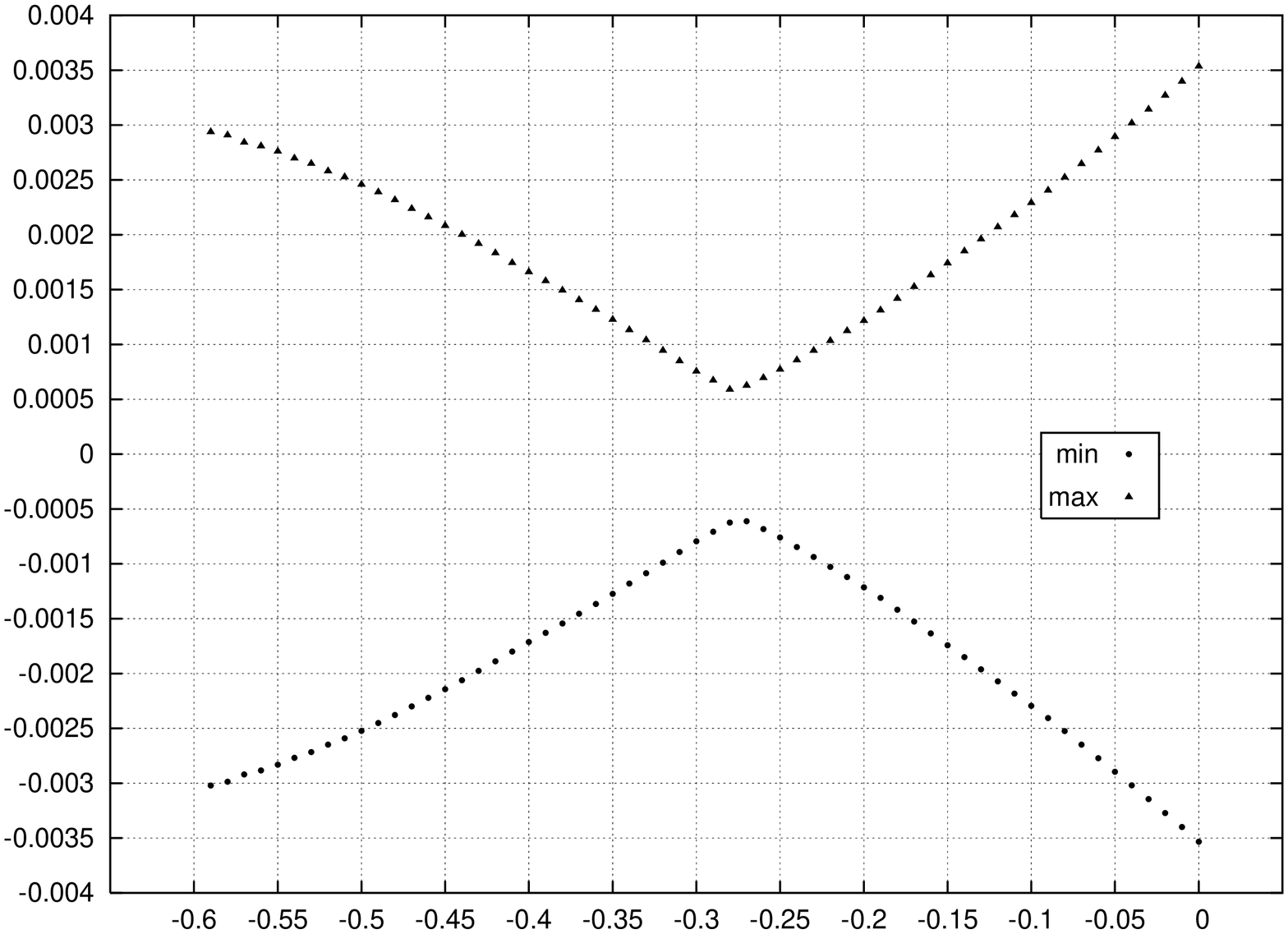} 
}
$$ \
\SetLabels
\L (0.48*0.900) {\small{$\psi_{\max}$}} \\
\L (0.48*0.130) {\small{$\psi_{\min}$}} \\
\L (0.98*-0.01) {\tiny{$\alpha$}} \\
\endSetLabels
\leavevmode 
& \hspace{7mm}
$$
\strut\AffixLabels{
\includegraphics[width=38mm,angle=270]{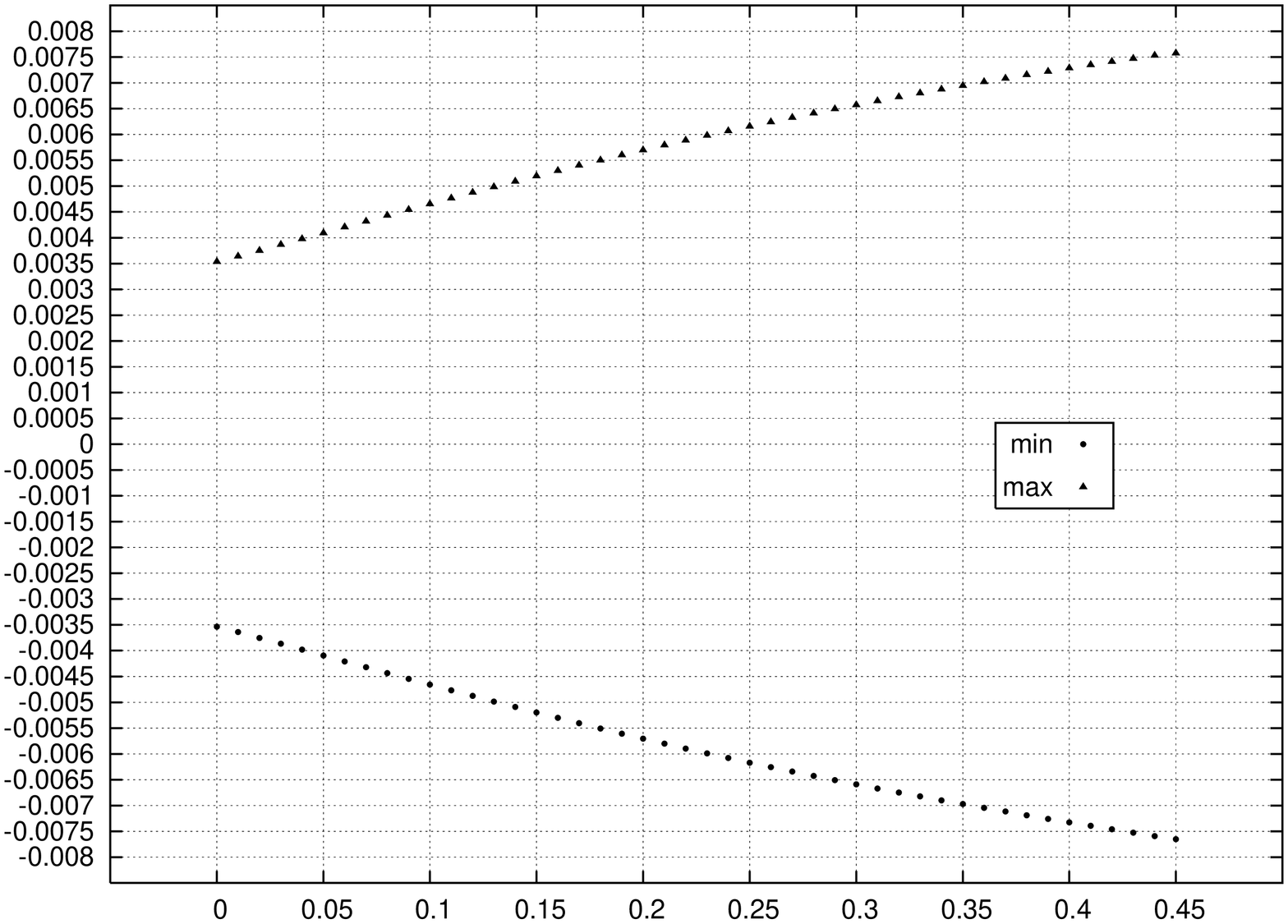} 
}
$$
 \vspace{1mm}\\
{\small (a)$\ \alpha$ negative} & 
{\small (b) $\ \alpha$ positive}
\end{tabular}
\caption{\small Profile of the extrema values of the stream function for inertial viscoelastic flows with 
${\cal R}e$=$5$ and $\delta$=$0.2$.}
\label{extrema_viscoinert}
\end{figure}

\begin{figure}[!ht]
\SetLabels
\L (0.48*0.99) $\rightarrow$ \\
\L (0.48*-0.01) $\rightarrow$ \\
\endSetLabels
\leavevmode
\centering
\begin{tabular}{ccc}
$$
\strut\AffixLabels{
\includegraphics[width=38mm]{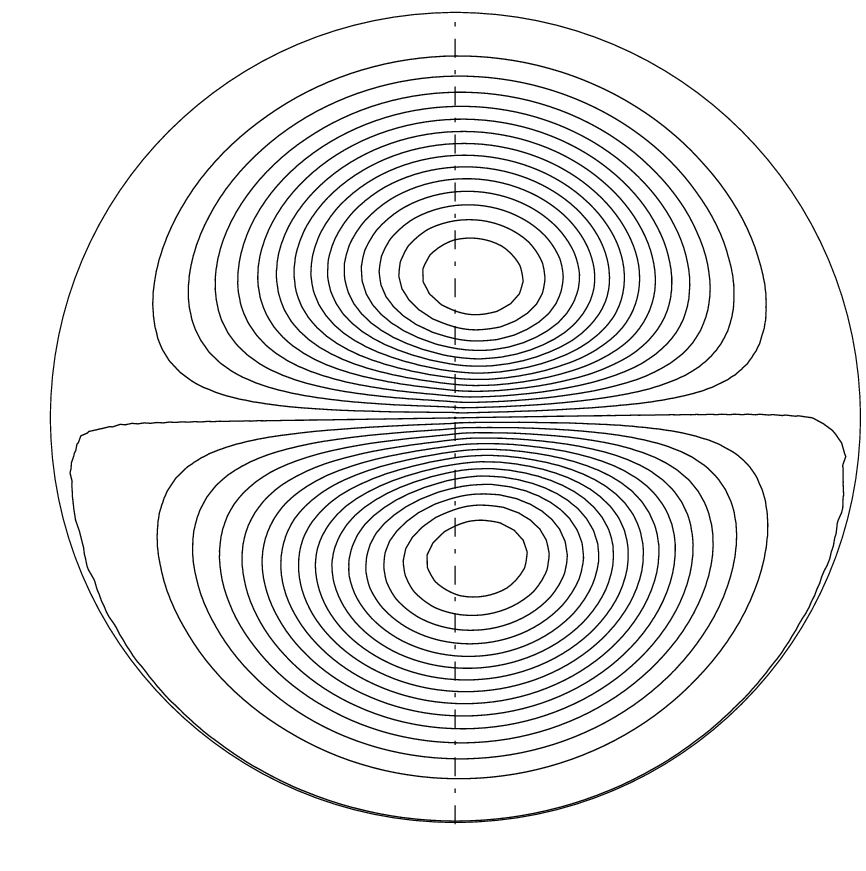} 
}
$$
\SetLabels
\L (0.48*0.99) $\leftarrow$ \\
\L (0.48*-0.01) $\leftarrow$ \\
\endSetLabels
\leavevmode
& 
\includegraphics[width=38mm]{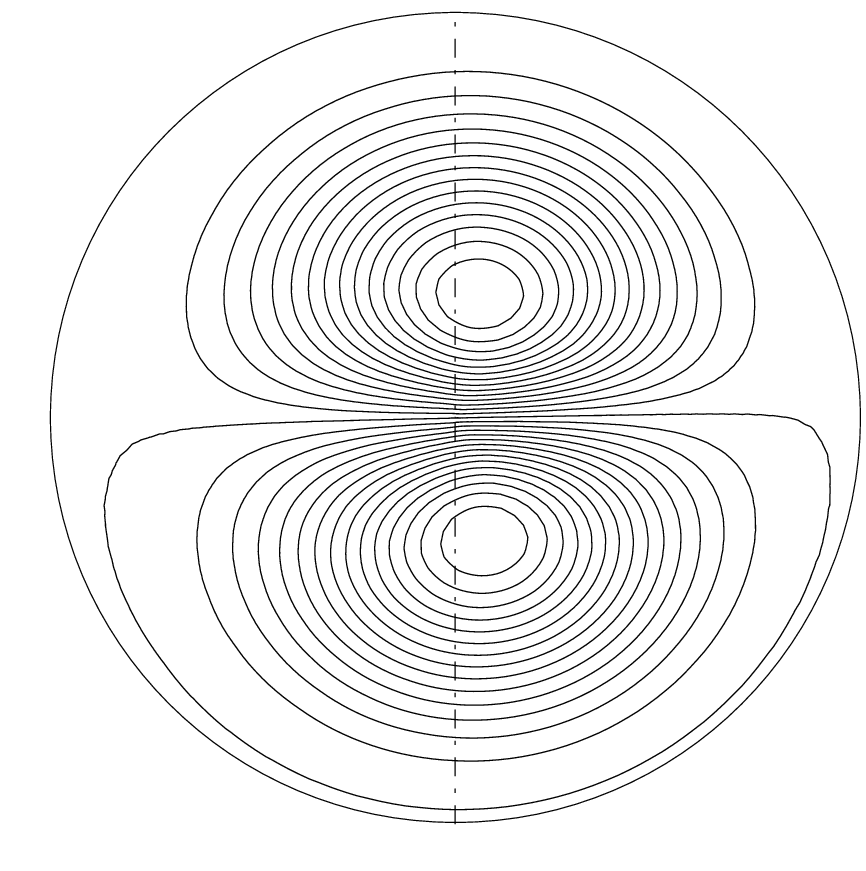} 
& \includegraphics[width=38mm]{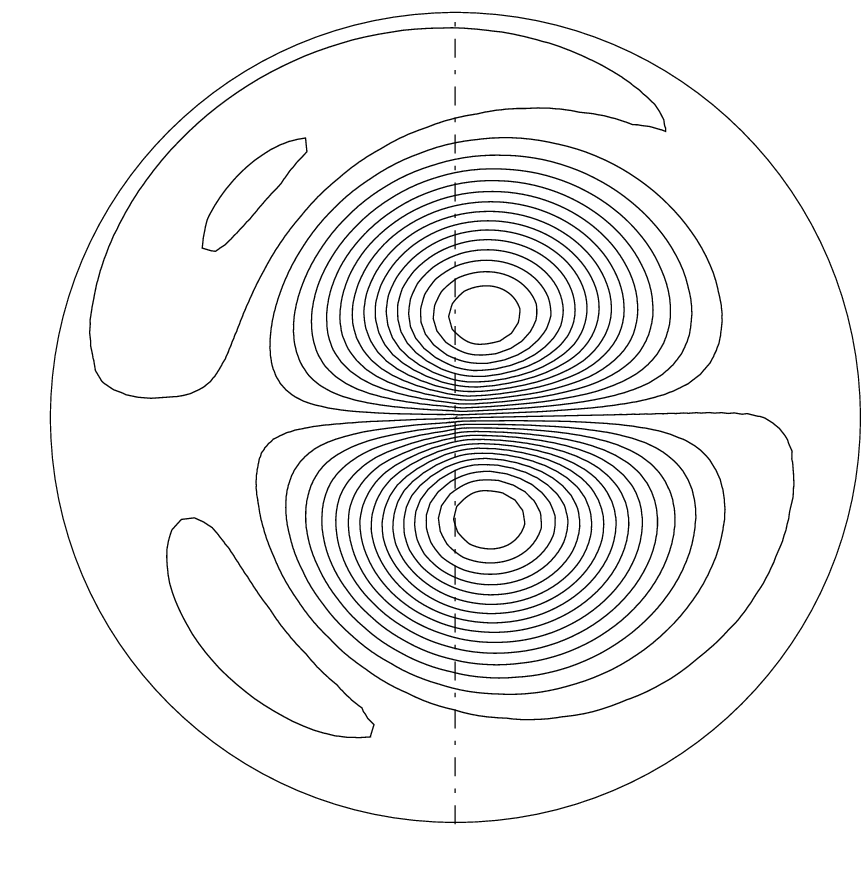}
\vspace{-2mm}\\
{\tiny{ (a) $\mathcal{R}e=0.5$}} & \hspace{4mm} {\tiny{ (b) $\mathcal{R}e=1.0$}} 
& \hspace{4mm} {\tiny{ (c) $\mathcal{R}e=1.4$}} \vspace{2mm}\\
\includegraphics[width=38mm]{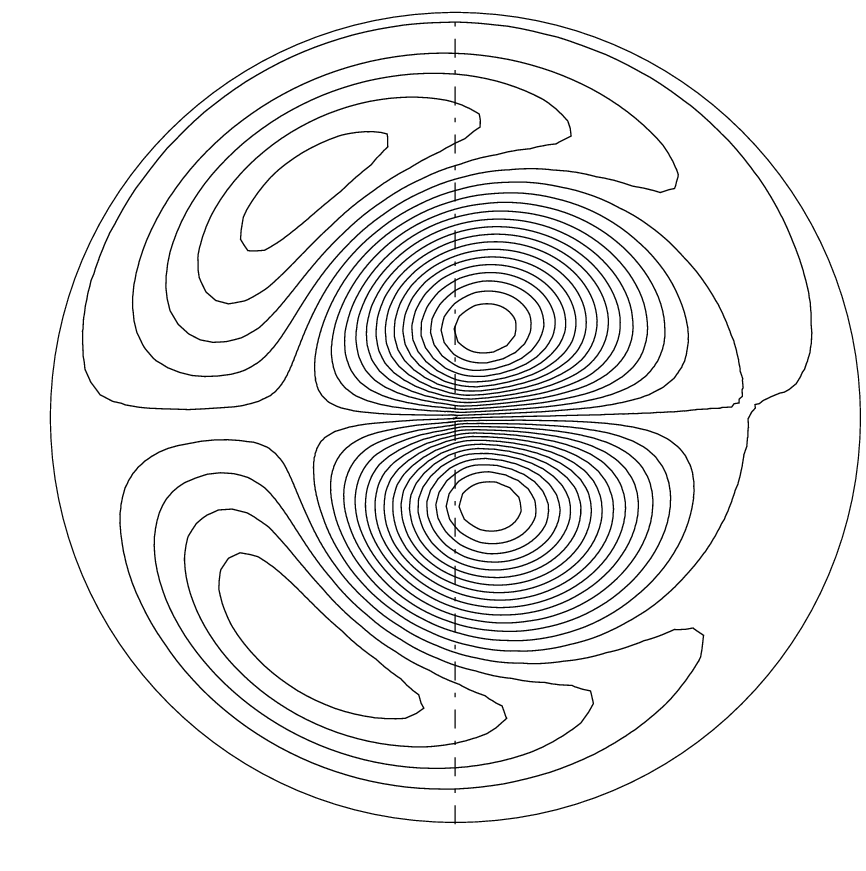} & 
\includegraphics[width=38mm]{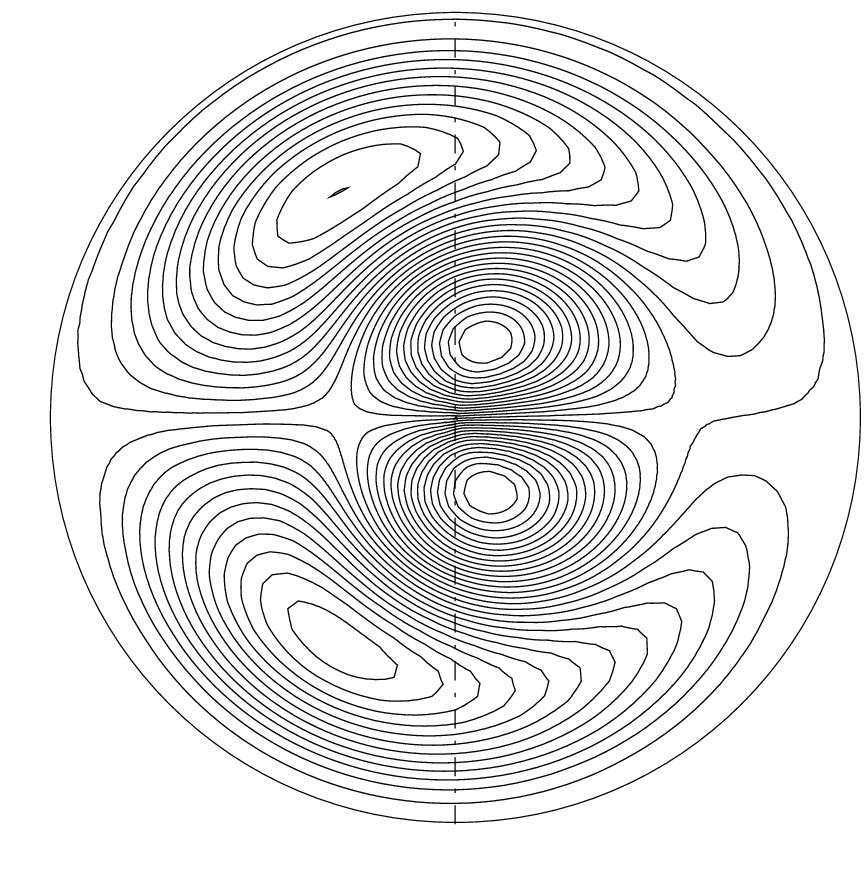}
& \includegraphics[width=38mm]{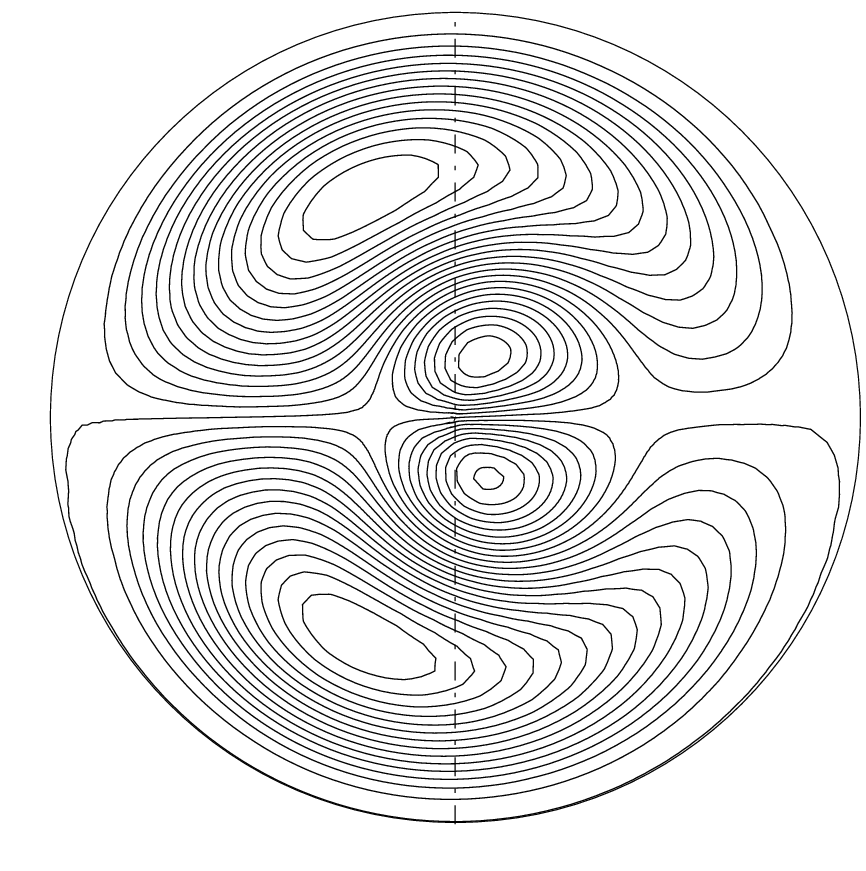}
\vspace{-2mm}\\ 
{\tiny{ (d) $\mathcal{R}e=1.6$}} & \hspace{4mm} {\tiny{ (e) $\mathcal{R}e=1.8$}}
& \hspace{4mm} {\tiny{ (f) $\mathcal{R}e=2.0$}} \vspace{2mm}\\
\includegraphics[width=38mm]{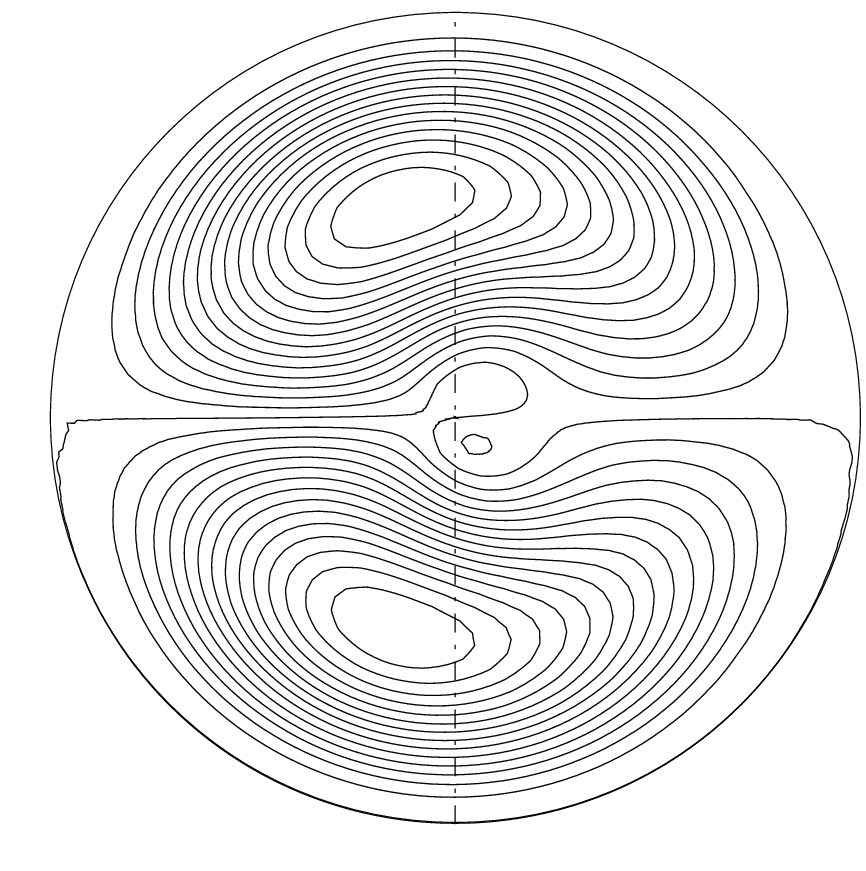} & 
\includegraphics[width=38mm]{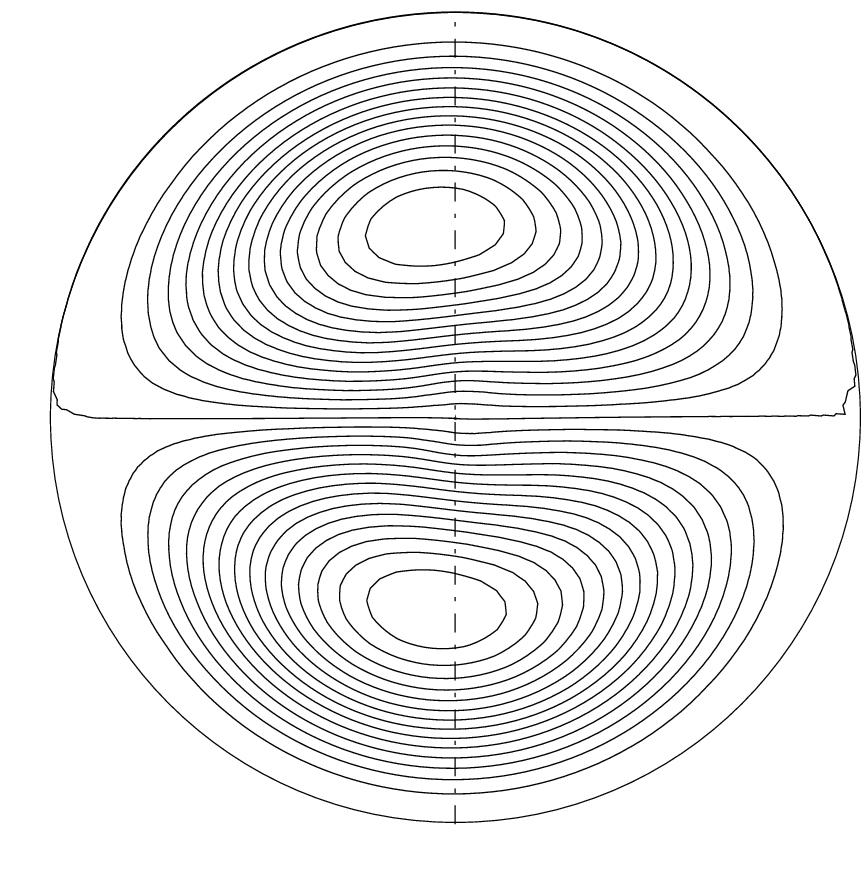} 
& 
$$
\strut\AffixLabels{
\includegraphics[width=38mm]{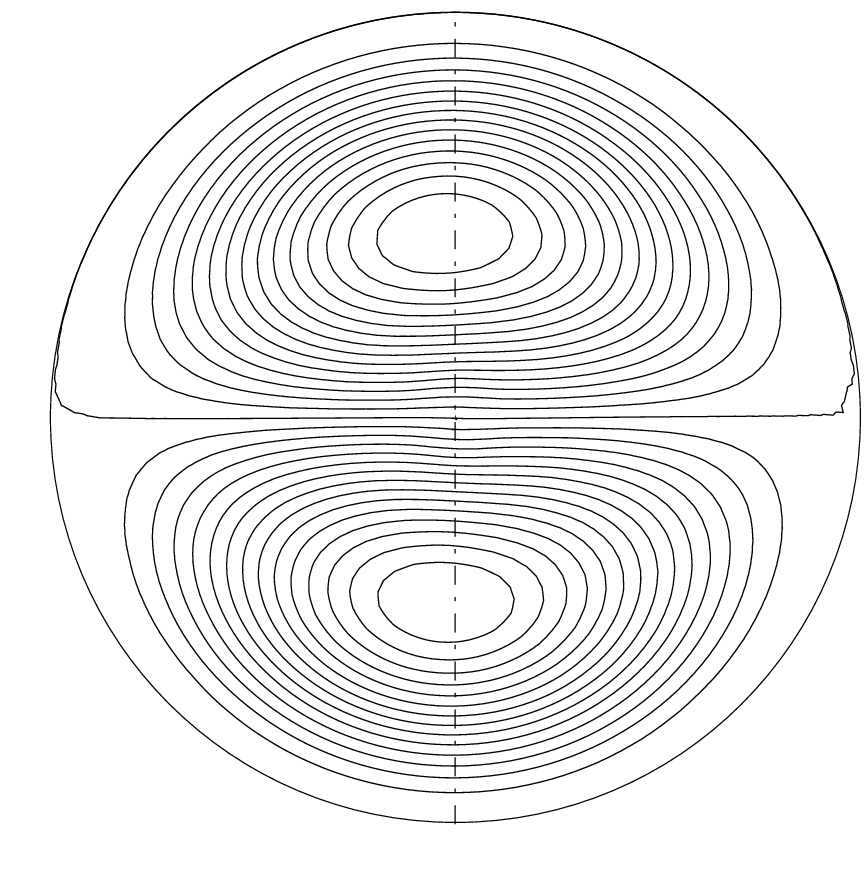}
}
$$ 
\vspace{-2mm}\\
{\tiny{ (g) $\mathcal{R}e=2.5$}} & \hspace{4mm} {\tiny{ (h) $\mathcal{R}e=5.0$}}
& \hspace{4mm} {\tiny{ (i) $\mathcal{R}e=10.0$}} \vspace{4mm}
\end{tabular}
\caption{\small Effects of the inertia on the streamlines with $\alpha=-0.1$ and  $\delta=0.2$.
}
\label{inert-visco0.5}
\end{figure}

\begin{figure}[!ht]
\SetLabels
\L (0.48*0.830) {\small{$\psi_{\max}$}} \\
\L (0.48*0.200) {\small{$\psi_{\min}$}} \\
\L (0.98*0.00) {\tiny{$\mathcal{R}\! e$}} \\
\endSetLabels
\leavevmode
\centering
\begin{tabular}{c}
$$
\strut\AffixLabels{
\includegraphics[width=40mm,angle=270]{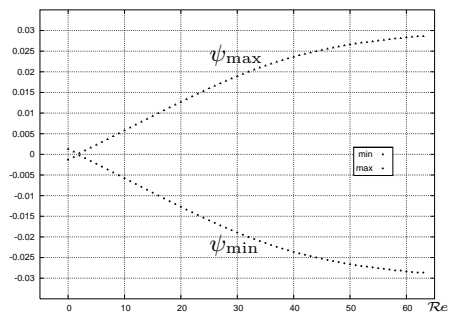}}
$$  
\end{tabular}
\caption{\small Profile of the extrema values of the stream function for inertial viscoelastic flows with 
$\alpha$=$-0.01$ and $\delta$=$0.2$.}
\label{extrema_inertvisco}
\end{figure}

\noindent Let us now consider the reversal phenomenon.
Fixing  the parameters $\alpha$, we study the effect
of fluid inertia.  Setting  $\alpha$=$-0.1$ and $\delta$=$0.2$, we
present in Figure \ref{inert-visco0.5} the contour plots of the stream function for different values of 
${\cal R}e$ in the interval $[0.5, 4]$. 
For ${\cal R}e=0.5$,
the viscoelasticity is dominant and the inertial forces have no real 
effect on the secondary flows. 
For ${\cal R}e=1.4$,
some modifications occur. The streamlines in the core region appear to be less dense,
 the size of the couple of 
vortices is smaller, the flow is driven near the wall pipe. Moreover, 
we observe the formation of boundary layer flows, with a pair of weak and elongated vortices.
As Reynolds number increases, the boundary layer flows increase rapidly.
The vortices in the core region are weaker, while 
the streamlines near the pipe wall are more dense and distorted, and
a strenghthening of the new vortices becomes clear. There is evidence that 
in contrast with the case where the viscoelasticity dominates,
 the wall pipe here 
is the transition region from the viscoelastic state to the inertial one,
 and that the reversal flows, formed near the
boundary around the new vortices, develop and drive the flow to the center
 of the pipe. 
Finally, we can conclude that for sufficiently large Reynolds number the secondary 
flows behave as in the Newtonian case superimposing the $\alpha$ parameter 
influence.
\section{Conclusion}
The finite element numerical simulations presented in this work provide
relevant information on the qualitative and quantitative behaviour of steady 
fully developed flows of second-grade fluids in curved pipes of circular cross 
section and arbitrary curvature ratio, driven by a pressure drop. As stated above, 
the main feature of curved pipe flows is the existence of secondary motions which 
are clearly related to small changes with respect to the different parameters of 
the problem. The numerical tests presented here show that secondary flows are 
promoted by both viscoelasticity and inertia with a globally newtonian behaviour 
when the fluids are compatible with the thermodynamics. The more complex case of 
negative viscoelastic parameters show that the viscoelasticity and the inertia 
have opposite effects.\vspace{5mm}\\
\sloppy
\noindent{\bf \large Acknowledgements.}  The second author was supported by Centro de Investiga\c{c}\~ao em Matem\'atica e Aplica\c{c}\~oes (CIMA) through the grant UID/MAT/04674/2013 of FCT-Funda\c{c}\~ao para a Ci\^encia e a Tecnologia.


\begin{thebibliography}{99}
\vspace{4mm}
\baselineskip 4mm

\bibitem {ABCH05}
M.~Amara, C.~Bernardi, V.~Girault and F.~Hecht, Regularized finite element discretizations of a 
grade-two fluid model, {Int. J. Num. Meth. Fluids}~48 (2005) 1375--1414.

\vspace{-3pt}
\bibitem {citation1} 
N.~Arada, P.~Correia and A.~Sequeira, Analysis and finite element simulations of a second-order 
fluid model in a bounded domain, {Numerical Methods for Partial Differential Equations}, 23 (2007) 1468--1500.

\vspace{-3pt}
\bibitem{arada-pires-sequeira1} N. Arada, M. Pires and A. Sequeira, Viscosity effects on flows of generalized 
Newtonian fluids through curved pipes, {Computers and Mathematics with 
Applications}, 53 (2007) 625--646.

\vspace{-3pt}
\bibitem{arada-pires-sequeira2} N. Arada, M. Pires and A. Sequeira, 
Numerical simulations of shear-thinning Oldroyd-B fluids in curved pipes,
IASME Transactions, Issue 6, 2 (2005) 948-959. 

\vspace{-3pt}
\bibitem{berger} S. A. Berger, L. Talbot and L.-S. Yao,  Flow in curved pipes,
{Ann. Rev. Fluid Mech.}, 15 (1983) 461-512.

\vspace{-3pt} 
\bibitem{BAH87} R. B. Bird, R. C. Armstrong and O. Hassager,
{Dynamics of polymeric liquids}, John Wiley \& Sons, New York (1987).
 
\vspace{-3pt}
\bibitem{barnes} H. A. Barnes and K. Walters,
On the flow of viscous and elastico-viscous liquids through straight
and curved pipes, { Proc. Roy. Soc. Lond.}, A 314 (1969) 85--109.

\vspace{-3pt}
\bibitem{bowen} 
P. J.~Bowen, A.R.~Davies and K.~Walters,
On viscoelastic effects in swirling flows,
{J. Non--Newtonian Fluid Mech.}~38 (1991) 113--126.

\vspace{-3pt}
\bibitem{correia} 
P.~Correia, {Numerical simulations of a non--Newtonian fluid
flow model using finite element methods}, 
PhD Thesis, IST, Lisbon, 2004.

\vspace{-3pt}
\bibitem{CoRo} 
V.~Coscia and A.~Robertson,  Existence and uniqueness of steady,
fully developed flows of second order fluids in curved pipes,
{Math. Models Methods Appl. Sci.}~11 (2001) 1055--1071.

\vspace{-3pt}
\bibitem{dean1} 
W.R.~Dean, Note on the motion of fluid in curved pipe,
{Philos. Mag.}~4 (1927) 208--223.

\vspace{-3pt}
\bibitem{dean2} 
W.R.~Dean, The streamline motion of fluid in curved pipe,
{Philos. Mag.}~5 (1928) 673--695.

\vspace{-3pt}
\bibitem{fan} Y. Fan, R. I. Tanner and N. Phan-Thien, Fully developed viscous and viscoelastic flows in curved pipes, 
{ J. Fluid Mech.}, 440 (2001) 327-357.

\vspace{-3pt}
\bibitem {citation2} 
G.P.~Galdi and A.~Sequeira, Further existence results for classical
solutions of the equations of second--grade fluid,
{Arch. Rational Mech. Anal.}~128 (1994) 297--312.

\vspace{-3pt}
\bibitem{GiSc02} V.~Girault and L.~R.~Scott, Finite element discretizations of a 
two-dimensional
grade-two fluid model, {M2AN}~35 (2002) 1007--1053.

\vspace{-3pt}
\bibitem{ito} H. Ito, Flow in curved pipes, { JSME Int. J.},
 30 (1987) 543-552.

\vspace{-3pt}
\bibitem {citation3} 
W.~Jitchote and A.M.~Robertson, Flow of second order fluids in curved pipes,
{J. Non--Newtonian Fluid Mech.}~90 (2000) 91--116.

\bibitem{Rob} A.~M.~Robertson, On viscous flow in curved pipes of non-uniform cross section,
{Inter. J. Numer. Meth. fluid.}~22 (1996) 771--798.

\vspace{-3pt}
\bibitem{sohr} 
W.Y.~Soh and S.A.~Berger, Fully developed flow in a curved pipe of
ar\-bitrary curvature ratio,
{Int. J. Numer. Meth. Fluid.}~7 (1987) 733--755.

\end{thebibliography}
\end{document}